%% file: main.tex
\documentclass{article}

\usepackage{arxiv}

\usepackage[utf8]{inputenc} 
\usepackage[T1]{fontenc}    
\usepackage{hyperref}       
\usepackage{url}            
\usepackage{booktabs}       
\usepackage{nicefrac}       
\usepackage{microtype}      
\usepackage{lipsum}
\usepackage{graphicx}
\usepackage{import}
\usepackage{float}
\usepackage{cite}
\usepackage{amsmath,amssymb,amsfonts}
\usepackage{algorithm,algpseudocode}
\usepackage{booktabs}
\usepackage{array}
\usepackage{lineno,hyperref}
\usepackage{textcomp}
\usepackage{relsize}
\usepackage{subcaption}
\graphicspath{ {./images/} }

\algblock{ParFor}{EndParFor}
\algnewcommand\algorithmicparfor{\textbf{parfor}}
\algnewcommand\algorithmicpardo{\textbf{do}}
\algnewcommand\algorithmicendparfor{\textbf{end\ parfor}}
\algrenewtext{ParFor}[1]{\algorithmicparfor\ #1\ \algorithmicpardo}
\algrenewtext{EndParFor}{\algorithmicendparfor}

\title{A Novel Highly Parallelizable Machine-Learning Based Method for the Fast Solution of Integral Equations for Electromagnetic Scattering Problems}

\author{
 Enes~Koç \\
  Department of Electrical and Electronics Engineering\\
   Bilkent University\\
  Ankara, 06800, Türkiye \\
  \texttt{enes.koc@bilkent.edu.tr} \\
   \And
 Mert~Kalfa \\
  Department of Electrical and Electronics Engineering\\
   Bilkent University\\
  Ankara, 06800, Türkiye \\
  \AND
 Secil~E.~Dogan \\
  Department of Electrical and Computer Engineering\\
  The Ohio State University\\
  Columbus, OH 43210, USA \\
  \And
   Vakur~B.~Ertürk \\
  Department of Electrical and Electronics Engineering\\
   Bilkent University\\
  Ankara, 06800, Türkiye \\
}

\begin{document}
\maketitle
\begin{abstract}
We propose a novel method for the efficient and accurate iterative solution of frequency domain integral equations (IEs) that are used for large/multi-scale electromagnetic scattering problems. The proposed method uses a novel group-by-group interaction strategy to accurately evaluate far-zone interactions within the framework of the one-box-buffer scheme during the matrix-vector multiplication at each iteration. Briefly, subdomain basis functions that are used to model the scatterer at each box are represented by a fixed number of uniformly distributed and arbitrarily oriented Hertzian dipoles (referred to as uniform basis functions), and then the dipole-to-dipole interactions are predicted in a group-wise manner by employing machine learning algorithms, thereby showcasing efficiency, strong scalability for parallelization and accuracy without the low-frequency breakdown (LFB) problem. Since the dipole representation is independent of the underlying material properties of the scatterer, the proposed method is valid for all types of IEs (surface or volume). Moreover, because the training is performed offline, the resulting networks can be used for any scatterer under any IE, without extra training, as long as the size of, and the distances among the boxes are preserved. The efficiency and accuracy of the proposed method are assessed by comparing our results with those obtained from the conventional multilevel fast multipole algorithm for various scattering problems. The proposed method's parallelization performance is showcased through scalability tests, and its resilience to LFB is demonstrated.
\end{abstract}


\section{Introduction}
\label{sec:introduction}
Integral equation (IE) solvers in frequency domain using method of moments (MoM)~\cite{harrington1996field}, or its accelerated versions such as the fast multipole algorithm (FMM)~\cite{FMM}, and the multilevel fast multipole algorithm  (MLFMA)~\cite{song1997multilevel} are among the most widely employed methods for addressing electromagnetic (EM) radiation and scattering problems. Notably, MLFMA stands out as one of the foremost enhancements by achieving a computational complexity of $\mathcal{O}(N\log{N})$ for $N$ unknowns (i.e., discretization elements). Growing requirements to investigate more complex, large/multi-scale scatterers have led to many efficiency and accuracy improvements on the conventional MLFMA.

Regarding the efficiency of MLFMA, apart from well-crafted novel tree structures~\cite{IL-MLFMA1,IL-MLFMA2} and compelling skeletonization techniques~\cite{skeleton}, sophisticated state-of-the-art parallelization schemes~\cite{parallel2,parallel} have vastly enhanced its performance. However, achieving effective load balancing in parallelization of MLFMA is problematic due to irregular computational workload, which varies with the geometry of the problem. Therefore, even the most advanced schemes lack strong scalability and require very careful implementations to be efficient. Recently, machine learning (ML) techniques~\cite{Goodfellow}, which have become a valuable tool in various areas of electromagnetics~\cite{antenna1,antenna2,antenna3,invscat1,invscat2,intpred}, have been used to enhance the efficiency of IE-based solvers, particularly in scattering problems~\cite{time1,time2,Halil,MLTranslation}.

Regarding the accuracy of the conventional MLFMA, various strategies have been proposed to mitigate the well-known low-frequency breakdown (LFB) problem~\cite{LFB1}. When the translation distance between the far-zone MLFMA boxes is electrically small, calculating the interactions among them using multipole expansions~\cite{LFB2,LFB3,LFB5} or deforming the angular integration path in order to capture the evanescent wave contributions~\cite{LFB1,LFB6,LFB7,LFB8} are some of the available remedies. Unfortunately, many of these remedies require the solver to be implemented from the ground up with increased complexity due to the necessity of an alternative expansion of the Green's function. An approximate diagonal form~\cite{OzgurBariscan} is used in the conventional MLFMA, which has a simple implementation. Unfortunately, it is not as accurate as the previously reported works. Alternatively, use of multiple precision arithmetic in MLFMA has been proposed in~\cite{LFB9} to mitigate the LFB problem with the downside of requiring specialized hardware implementations.

On the other hand, growing requirements to investigate more complex, large/multi-scale scatterers also underscore the necessity of a search for novel computational methods with increased efficiency. Therefore, in this paper we propose a novel, highly parallelizable and broadband method for the fast and accurate iterative solution of frequency domain IEs that has a potential to circumvent the aforementioned problems. The backbone of the proposed method is a novel group-by-group interaction strategy for the fast and accurate evaluation of far-zone interactions using the one-box-buffer scheme during the matrix-vector multiplications (MVMs) at each iteration. Briefly, subdomain basis functions that are used to discretize any scatterer at each and every identical box are represented by a fixed number of arbitrarily oriented and uniformly distributed Hertzian dipoles, referred to as uniform basis functions\footnote{The different choices of uniform basis functions are allowed in the formulation; however, Hertzian dipoles are selected due to their simple Green's function representations.}. Then, dipole-to-dipole interactions are calculated in a group-wise manner by employing an ML algorithm that predicts the interactions by using the Hertzian dipole's Green's function (HDGF) in its training stage. Consequently, \textit{i)} because of the workload for box-to-box interactions is now identical and independent throughout the entire geometry, a strong scalability for parallelization is achieved; \textit{ii)} the use of HDGF during the training of the ML networks avoids the LFB problem leading to an accurate broadband solution that can be employed to multi-scale scatterers. Moreover, because the training is performed offline, the resulting networks can be used for any scatterer under any surface or volume IE, without extra training, as long as the sizes of, and the distances among the boxes are preserved. The efficiency and accuracy of the proposed method are validated for both surface and volume IEs by comparing our scattering results for various objects with those obtained from the conventional MLFMA (and Mie series when applicable), and its resilience to LFB is confirmed. Furthermore, scalability tests demonstrate the parallelization performance of the proposed method.

The remainder of this article is organized as follows: Section~\ref{sec2} discusses the formulation with a focus on the novel group-by-group interaction model, which forms the foundation for the proposed method. Section~\ref{sec3} explains the ML models used to characterize the group-by-group interactions. Next, the parallel implementation scheme of the proposed method is given in Section~\ref{sec4}. Numerical results are presented in Section~\ref{sec5}. Finally, Section~\ref{sec6} provides concluding remarks. An $\exp{(j\omega t)}$ time convention, where $\omega = 2\pi f$ and $f$ representing the operating frequency, is adopted and suppressed throughout the paper.   

\section{Formulation} \label{sec2}
In a typical MoM formulation, an IE (volume or surface) can be cast into a matrix equation in the form of 
\begin{equation}
\label{eqn:1}
    \bar{\bar{Z}} \bar{I}=\bar{V},
\end{equation}
where $\bar{\bar{Z}}$ is the known impedance matrix, whose entries represent the interactions between basis \& test functions that are identical with a Galerkin scheme, $\bar{V}$ is the known excitation vector, and finally $\bar{I}$ contains the unknown coefficients of the subdomain basis functions to be found iteratively by employing the novel method proposed in this paper. 

\begin{figure}[H] 
\centering
\includegraphics[width=3.4in]{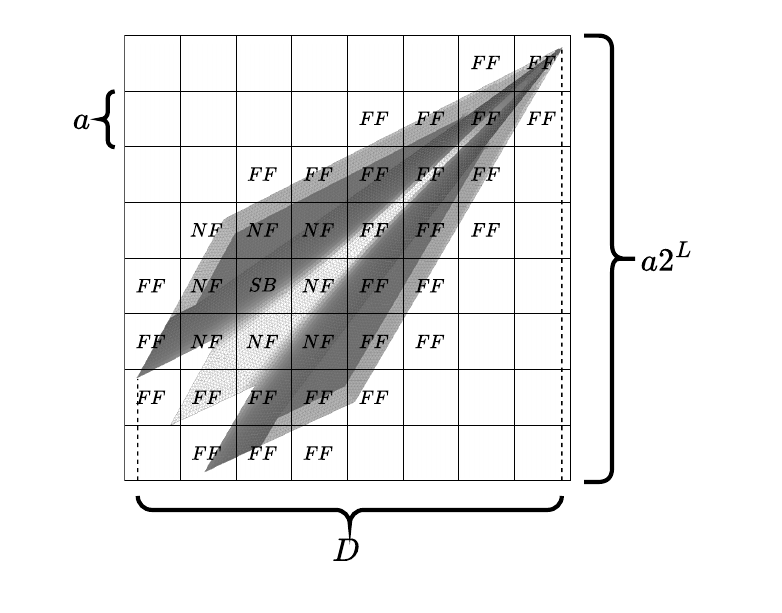}%
\caption{2D projection of a sample scatterer inside a mother box, which is divided into child boxes. For the selected box $SB$, $NF$ boxes correspond to near interactions, and $FF$ boxes correspond to far interactions.}
\label{fig1}
\end{figure}

Similar to FMM, our proposed method starts by enclosing the entire object, sized $D$, with a mother box (cube) and dividing it into $8^{L}$ ($L:$ a positive integer) equal sized child boxes each of which has a side length of $a$. Hence, the edge length of the mother box is $a2^{L}$, and the inequality $a2^{L-1}<D\leq a2^{L}$ should hold for all objects. Then, using the one-box-buffer scheme, as shown in Fig.~\ref{fig1}, and only considering the non-empty boxes, MVM operations (at each iteration) can be separated into near (N) and far (F) interactions given by 
\begin{equation}
\label{eqn:2}
    \bar{\bar{Z}}^{N} \bar{I} + \bar{\bar{Z}}^{F} \bar{I} = \bar{V}.
\end{equation}
In~\eqref{eqn:2}, $\bar{\bar{Z}}^{N}$ denotes the near-interaction matrix that is calculated directly and stored in the memory, while $\bar{\bar{Z}}^{F}$ stands for the far-zone interaction matrix and is calculated by using the proposed novel group-by-group interaction computation strategy.  The group-by-group interaction scheme has three distinct stages, namely \textit{1)} mapping to uniform basis functions, \textit{2)} translation, and \textit{3)} inverse mapping to test (subdomain) functions, which are explained in the following subsections.  

\subsection{Mapping to Uniform Basis Functions}
\begin{figure}[H]
\centering
\includegraphics[width=4.5in]{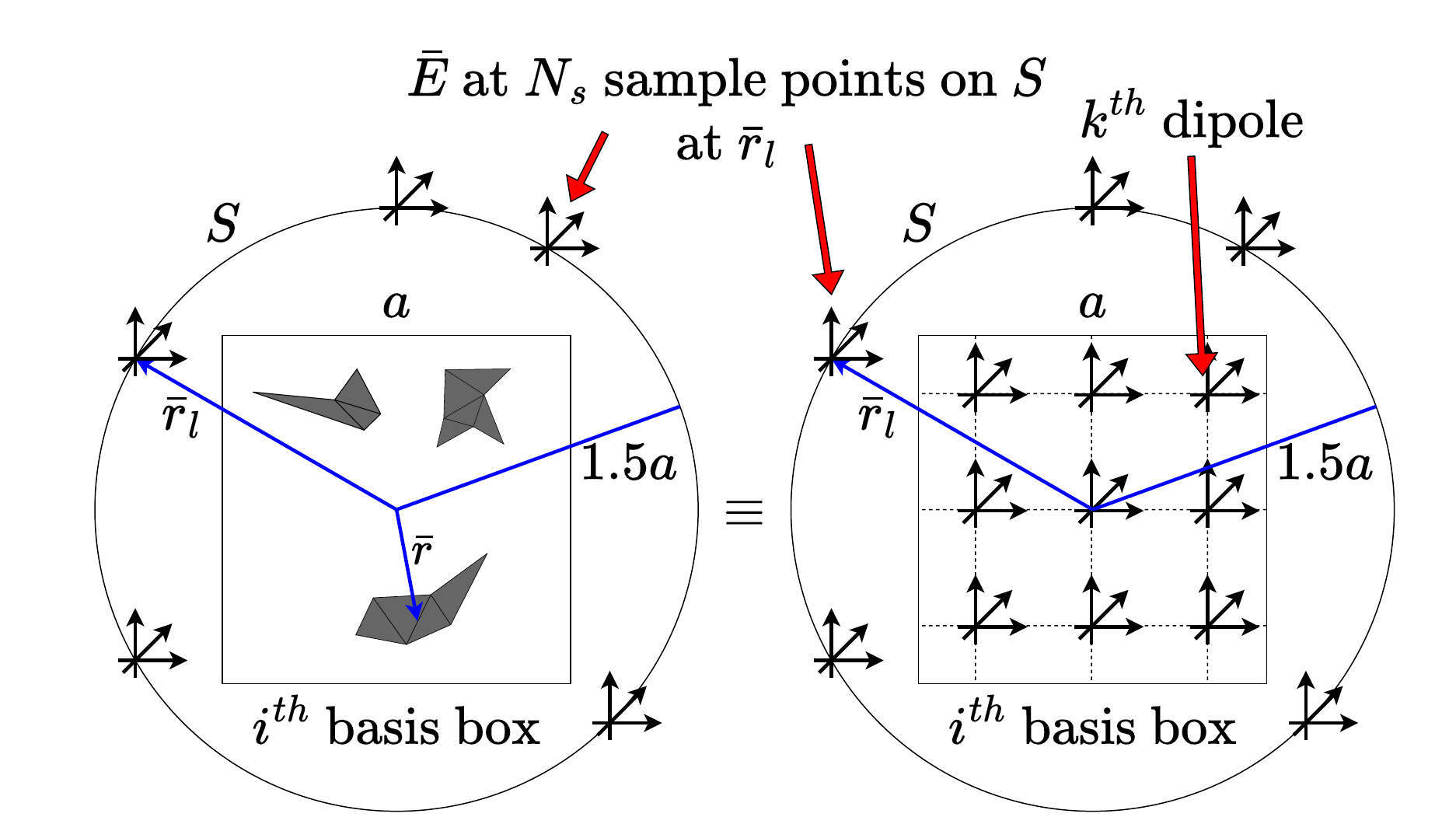}
\caption{Mapping from subdomain basis functions to uniform basis functions} 
\label{fig2}
\end{figure}

The goal of this stage is to represent the subdomain basis functions that are used to discretize an arbitrary scatterer with any material property in each box with a fixed number of arbitrarily oriented and uniformly distributed basis functions that we call uniform basis functions. In this paper, arbitrarily oriented Hertzian dipoles (i.e., each dipole is a linear combination of $x-$, $y-$, and $z-$ directed dipoles) are used as uniform basis functions, as illustrated in Fig.~\ref{fig2}. 

Let $\bar{J}(\bar{r})$ represent the current distribution inside one of the child boxes given by
\begin{equation}
\label{eqn:3}
    \bar{J}(\bar{r}) = \sum_{k=1}^{N_b} a_{k}\bar{b}_k(\bar{r}),
\end{equation}
where $\bar{b}_k(\bar{r})$ is the original subdomain basis function with $a_k$ being its known coefficient (to be updated at each iteration) and $N_b$ is the number of subdomain basis functions inside the selected box. This current distribution creates the original electric fields given by
\begin{equation}
\label{eqn:4}
    \bar{E}_o(\bar{r}) = \int \bar{\bar{D}}(\bar{r},\bar{r}^{\prime}) \cdot \bar{J}(\bar{r}^{\prime}) d\bar{r}^{\prime},
\end{equation}
where $\bar{\bar{D}}(\bar{r},\bar{r}^{\prime})$ is the free-space Dyadic Green's function expressed as 
\begin{equation}
\label{eqn:5}
\begin{gathered}
    \bar{\bar{D}}(\bar{r},\bar{r}^{\prime}) = \psi(R)\bigg\{\bar{\bar{I}}_{3\times3}\left(1-\frac{j}{k R}-\frac{1}{k^2 R^2}\right)\\ - \hat{R} \hat{R}\left(1-\frac{3 j}{k R}-\frac{3}{k^2 R^2}\right)\bigg\}.
\end{gathered}
\end{equation}
In~\eqref{eqn:5}
\begin{equation}
    \label{eqn:6}
    \psi(R) = -j k \eta \frac{\exp (-j k R)}{4 \pi R},\\
\end{equation}
with $\bar{R} = \bar{r}-\bar{r}^{\prime},\ R = \|R\|,\ \hat{R} = \bar{R}/R,$ and $k$ and $\eta$ are the free-space wave number and intrinsic impedance, respectively. Similarly, let $\bar{J}_{u}(\bar{r})$ represent the current distribution of the uniform basis functions given as 
\begin{equation}
\label{eqn:7}
    \bar{J}_u(\bar{r}) = \sum_{k=1}^{N_u} c_{k}\bar{p}_{d,k}(\bar{r}),
\end{equation}
where $\bar{p}_{d,k}$ is the position fixed uniform basis function, $c_k$ is its unknown coefficient to be determined at the end of this stage at each iteration, and $N_{u}$ is the number of uniform basis functions in each and every box (that is also fixed). The fields created by $\bar{J}_u(\bar{r})$ can also be expressed as 
\begin{equation}
    \label{eqn:8}
    \bar{E}_u(\bar{r}) = \sum_{k=1}^{N_u} c_{k}\int \bar{\bar{D}}(\bar{r},\bar{r}^{\prime}) \cdot \bar{p}_{d,k}(\bar{r}^{\prime}) d\bar{r}^{\prime}.
\end{equation}
The uniqueness theorem suggests that if two $\bar{E}$-fields over a closed surface $S$ are the same, with each of them created by the sources inside of $S$, then $\bar{E}$-fields outside of $S$ are the same. Therefore, we seek to find the best coefficients for the uniform basis functions (i.e., $c_k:k=1,\ldots,N_{u}$) by minimizing the $\bar{E}$-field error on the selected surface $S$, which is chosen as a sphere with a radius $1.5a$ as shown in Fig.~\ref{fig2}. The rationale for this choice is that we want each far-zone interaction to be an exterior point to $S$, and we select the biggest sphere that does not cross into test boxes within the one-box buffer scheme. As a result, our objective is to find $c_{k}$'s such that
\begin{equation}
\label{eqn:9}
    \underset{c_{k}}{\mathrm{argmin}} \int_{\bar{r}\in S} \| \bar{E}_o(\bar{r}) - \bar{E}_u(\bar{r})\|^2\ d\bar{r},
\end{equation}
which can be discretized by choosing $N_s$ sampling points on the sphere $S$, and~\eqref{eqn:9} can be rewritten as 
\begin{equation}
\label{eqn:10}
    \underset{c_{k}}{\mathrm{argmin}} \sum_{l=1}^{N_{s}} \| \bar{E}_o(\bar{r}_l) - \bar{E}_u(\bar{r}_l)\|^2 = \underset{c_{k}}{\mathrm{argmin}}  \|\bar{B} - \bar{\bar{A}} \bar{C}\|^{2}.
\end{equation}
In~\eqref{eqn:10}, $\bar{B}$ is $3N_s \times 1$ vector that contains components of $\bar{E}_o$ at sample points (i.e., $\bar{r}_{l}\in S$), $\bar{\bar{A}}$ is the $3N_s \times N_u$ matrix, whose entries represent the $\bar{E}$-field contribution of $\bar{p}_{d,k}$'s at sample points\footnote{A factor 3 appears in the sizes of $\bar{\bar{A}}$ and $\bar{\bar{B}}$ due to independent $x,y,z$ components of the dipoles.} and finally $\bar{C}$ is the $N_u \times 1$ unknown coefficient vector of $c_k$'s. The least-squares solution of~\eqref{eqn:10} is obtained from 
\begin{equation}
\label{eqn:11}
    \bar{C} = \bar{\bar{A}}^{\dagger}\bar{B}
\end{equation}
where $\bar{\bar{A}}^{\dagger}$ is the pseudo-inverse of $\bar{\bar{A}}$.

In each child box, a total of $n_{d}$ Hertzian dipoles is placed along each $x,y,z$ directions with uniform separation. The overall dipole count in each box is $N_d = n_d^{3}$, and given that each dipole has three free parameters (one for each Cartesian component), $N_u = 3N_d$. In the selection of $N_s$ and $N_u$, we are inspired by MLFMA and depending on the box size $a$ (or $ka$) we initially selected $N_d$ as $N_d = 2(1.73 k a+2.16\left(d_0\right)^{2 / 3}\left(k a\right)^{1 / 3})^2$~\cite{EBF} or $N_d = 2(14.14 d_0-7.17)^2$~\cite{EBFLFB}, where $d_0$ corresponds to the desired digits of accuracy in the Green's function. Then $N_s$ is selected to be twice of $N_d$, which is $N_s = 2N_d$. Finally, regarding the sampling scheme on $S$, spherical Fibonacci Lattice Points~\cite{Gonz_lez_2009} are used due to their approximate uniformity. In conventional spherical coordinates $(r,\theta,\phi)$, the Fibonacci lattice is defined for points $r_{l}=(1.5a,\theta_l, \phi_l)$ with $\theta_l = 2 \pi l / \varphi$, $\phi_l = \arccos{(1 - (2l+1)/N_s)}$, where $\varphi$ is the golden ratio.\footnote{The golden ratio is defined as $\varphi=0.5+0.5\sqrt{5}$.} The pseudocode shown in Algorithm~\ref{alg:cap1} is used to perform this mapping stage at each iteration when solving~\eqref{eqn:2} iteratively. It should be noted that $\bar{\bar{A}}^{\dagger}$ in~\eqref{eqn:11} and $\bar{B}_{k}$'s in~\eqref{eqn:11} for each subdomain basis function $\bar{b}_{k}$ are calculated and stored in the preprocessing stage, and are not included in the iterative solver computations. 
\begin{algorithm}[H]
\caption{Mapping to Uniform Basis Functions}\label{alg:cap1}
\begin{algorithmic}
  \For{$i \in$ Non-empty Child Box Indices}
    \State{Initialize $\bar{C}_{i} = 0$ for the $i$\textsuperscript{th} box}
    \For{$k \in$ Subdomain Functions in Box $i$}
        \State{Calculate~\eqref{eqn:11} using $a_{k}$'s given by the}
        \State{iterative solver, to obtain $\bar{C}_{k}^{\prime}$.}
        \State{$\bar{C}_{i} = \bar{C}_{i} + \bar{C}_{k}^{\prime}$.}
    \EndFor
  \EndFor
  \State \textbf{return} {$\bar{C}_{i}$'s}
\end{algorithmic}
\end{algorithm}

\subsection{Translation}
\begin{figure}[H]
\centering
\includegraphics[width=4.5in]{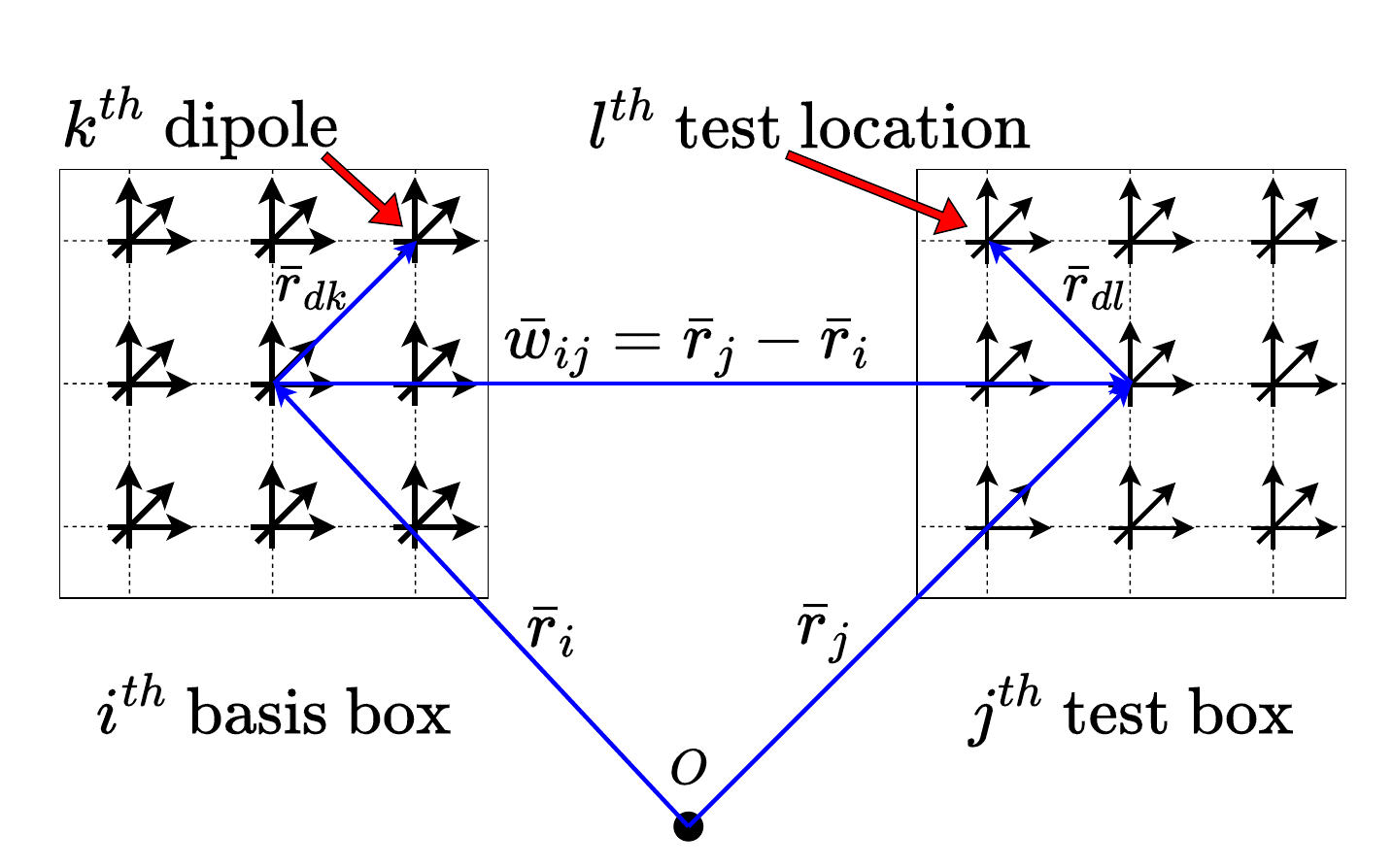}
\caption{Illustration of the translation stage.}
\label{fig3}
\end{figure}

Once the mapping to uniform basis functions stage is completed, all boxes become structurally identical (i.e., each box has the same number of fixed-position Hertzian dipoles).
Consider two boxes that are in the far-zone of each other, as shown in Fig.~\ref{fig3}, where $\bar{r}_{i}$, and $\bar{r}_{j}$ are the centers of the $i^{th}$ basis and $j^{th}$ test boxes, respectively, and $\bar{w}_{ij} = \bar{r}_{j}-\bar{r}_{i}$ is the translation vector from basis box to test box. Note that the translation vectors $\bar{w}$ are in the form of $a (\hat{x}n_1 +\hat{y}n_2+\hat{z}n_3)$ (shortened as $[n1,n2,n3]a$), where $-2^{L} < n_1,n_2,n_3 < 2^{L}$ and $n_1,n_2,n_3$ are integers. In other words, the child boxes are created to form a uniform lattice structure across the scatterer.  
Also in Fig.~\ref{fig3}, let $\bar{r}_{dk}$ be the $k^{th}$ dipole location with respect to the center of its box that is the $i^{th}$ basis box while let $\bar{r}_{dl}$ be the $l^{th}$ dipole location with respect to the center of its box that is the $j^{th}$ test box. The field at the $l^{th}$ dipole location (which is the test location) in the test box due to the $k^{th}$ dipole in the basis box is given by
\begin{equation}
\label{eqn:12}
    \bar{E}(\bar{r}_{dl}) = \bar{\bar{D}}(\bar{w}_{ij}+\bar{r}_{dl},\bar{r}_{dk})\cdot \bar{p}_{dk}.
\end{equation}
Generalizing~\eqref{eqn:12} for all dipoles and test locations (i.e., dipole locations) in the $i^{th}$ basis and $j^{th}$ test boxes, respectively, a matrix equation in the form of 

\begin{equation}
\label{eqn:13}
    \bar{E}_{i-j}(\bar{r}_{dl}) = \bar{\bar{G}}_{\bar{w}_{ij}}(\bar{w}_{ij}+\bar{r}_{dl},\bar{r}_{dk})\cdot \bar{P}_{i}(\bar{r}_{dk})
\end{equation}
can be obtained. In~\eqref{eqn:13}, $\bar{\bar{G}}_{\bar{w}_{ij}}$ is a $3N_d\times 3N_d$ translation matrix of the translation vector $\bar{w}_{ij}$. Each column of $\bar{\bar{G}}_{\bar{w}_{ij}}$ corresponds to the $\bar{E}$-field at test locations in the $j^{th}$ test box due to one component ($x$,$y$ or $z$) of a dipole inside the $i^{th}$ basis box without taking dipole's excitation coefficients into account (Hertzian dipole's Green's function). Also in~\eqref{eqn:13} $\bar{P}_{i}$ is the $3N_d\times 1$ vector corresponding to $\bar{C}$ that is solved via~\eqref{eqn:11}, and $\bar{E}_{i-j}$ is the $3N_d\times 1$ vector representing the fields at test locations. The total far-zone fields at test locations, denoted as $\bar{E}_{j}$, for $j^{th}$ test box is given as
\begin{equation}
    \label{eqn:14}
    \bar{E}_{j} = \sum_{i \in \mathbb{F}_{j}}\bar{E}_{i-j},
\end{equation}
where $\mathbb{F}_{j}$ is the set of the index of boxes that are in the far-zone of $j^{th}$ test box. At first glance, translation matrix $\bar{\bar{G}}_{\bar{w}}$ seems to be calculated for each translation vector $\bar{w}$ independently. However, since the boxes are structurally identical, many of the translation matrices can be derived from each other using rotational and mirror symmetries, provided that $|\bar{w}|$ is the same for the corresponding translation vectors. For example, if ${L=3}$, there are 316 distinct translation vectors, but only 16 unique translation matrices must be separately calculated since the remaining ones can be obtained via symmetries. 

The process of finding far-zone interactions by making use of $\bar{\bar{G}}_{\bar{w}}$ matrices is called the direct Green's function approach (DGFA) in this paper, and it is very accurate as it does not use the truncated multipole expansion or the truncated summation used in MLFMA. However, using DGFA in an element-by-element approach leads to a complexity of $\mathcal{O}(N_{d}^{2})$ for a single box-to-box interaction. Therefore, we employ ML networks that can accurately emulate translation matrices, thereby predicting all dipole-to-dipole interactions between any two far-zone boxes in a group-wise manner, which improves the efficiency significantly. Specific implementation of the ML networks is provided in Section~\ref{sec3}. Algorithm~\ref{alg:cap2} shows the pseudocode that performs the entire translation stage at each iteration using DGFA.
\begin{algorithm}[H]
\caption{Translation using DGFA}\label{alg:cap2}
\begin{algorithmic}
\State{Initialize $\bar{E}_{j}=0$ for $j=1,\ldots,N$}
  \For{$i \in$ Non-empty Child Box Indices}
    \For{$j \in$ Non-empty Child Box Indices}
        \If{$\|\bar{r}_{j}-\bar{r}_{i}\|=\|\bar{w}_{ij}\| \geq 2a$}
            \State{Calculate $\bar{E}_{i-j}$ defined in~\eqref{eqn:13}}
            \State{$\bar{E}_{j} = \bar{E}_{j} + \bar{E}_{i-j}$.}
        \EndIf
    \EndFor
  \EndFor
  \State \textbf{return} {$\bar{E}_{j}$'s}
\end{algorithmic}
\end{algorithm}
Note that the proposed ML networks will replace the DGFA calculation in~\eqref{eqn:13}.

\subsection{Inverse Mapping to Subdomain Functions}
At the end of the translation stage, fields at any test location can be determined via~\eqref{eqn:14}. As a matter of fact, these test locations can be selected as the interaction point locations of the original subdomain test functions. However, because we desire a scheme that is reusable, which does not depend on the object, and has an identical computational workload for each and every box from the parallelization point of view, we pursue a strategy similar to the first stage. That is, the field at any point inside of the test box is found from known fields at test locations by using the following mapping procedure. 

\begin{figure}[H]
\centering
\includegraphics[width=4.5in]{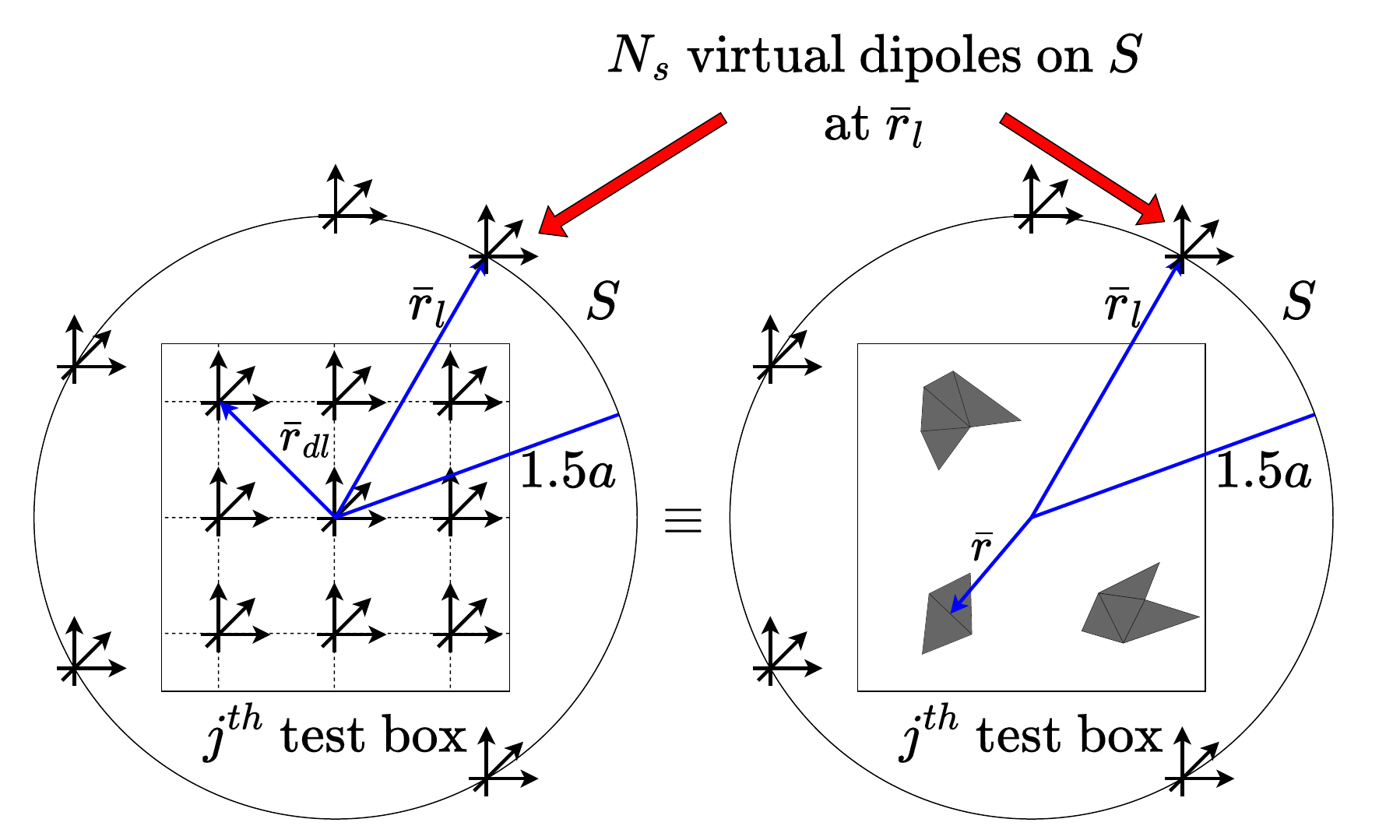}
\caption{Inverse mapping from uniform basis functions to subdomain functions}
\label{fig4}
\end{figure}

Let us consider the $j^{th}$ test box and place $N_{s}$ Hertzian dipoles (and call them virtual dipoles) on $S$ with $1.5a$ radius at the spherical Fibonacci lattice points, previously defined by $\bar{r}_{l}$ as shown in Fig.~\ref{fig4}.

For $j^{th}$ test box, fields inside at test locations $\bar{r}_{dl}$ generated by its far-zone boxes are given by 
\begin{equation}
\label{eqn:15}
    \bar{E}_{vd}(\bar{r}_{dl}) = \sum_{k=1}^{N_s} \bar{\bar{D}}(\bar{r}_{dl},\bar{r}_{l})\cdot q_{l}\bar{p}_{v,l}(\bar{r}_{l}),
\end{equation}
where $\bar{p}_{v,l}$ is the virtual dipole on $S$ at $\bar{r}_{l}$ with $q_{l}$ being its unknown excitation coefficient to be found. On the other hand, we already know the fields at these locations from the translation stage given by~\eqref{eqn:14} as $\bar{E_{j}}(\bar{r}_{dl})$. Thus, the problem of finding optimum $q_{l}$'s turns into another least-squares solution problem similar to that of the first stage, which can be formulated as 
\begin{equation}
\label{eqn:16}
\underset{q_{k}}{\mathrm{argmin}} \sum_{l=1}^{N_d} |\bar{E}_{j}(\bar{r}_{dl}) - \bar{E}_{vd}(\bar{r}_{dl}) |^{2} = \underset{q_{k}}{\mathrm{argmin}} \| \bar{E}_{j} - \bar{\bar{A}}^{\prime} \bar{Q}\|^{2},
\end{equation}
where $\bar{\bar{A}}^{\prime} = \bar{\bar{A}}^{T}$ is the same matrix used in~\eqref{eqn:10}, and $\bar{Q}$ is the $3N_s\times 1$ unknown coefficient vector of virtual dipoles. The least-squares solution of $\bar{Q}$ is obtained from 
\begin{equation}
    \label{eqn:17}
    \bar{Q} = (\bar{\bar{A}}^{T})^{\dagger}\bar{E}_{j}.
\end{equation}
Consequently, knowing $\bar{Q}$ and using~\eqref{eqn:15}, fields can be found anywhere inside the $j^{th}$ test box, including the desired integration points on the original subdomain functions via~\eqref{eqn:15}, except $\bar{r}_{dl}$ is now $\bar{r}$, as shown on the right side of Fig.~\ref{fig4}. Finally, the far-zone interaction matrix, and hence the MVMs are computed by testing the result of~\eqref{eqn:15} on the original subdomain functions. The pseudocode given in Algorithm~\ref{alg:cap3} provides the implementation of the inverse mapping stage at each iteration when solving~\eqref{eqn:2} iteratively. It should be noted that $(\bar{\bar{A}}^{T})^{\dagger}$ in~\eqref{eqn:11} and $\bar{\bar{D}}_{k}$'s, which are the vectorized versions of~\eqref{eqn:15} for $k^{th}$ subdomain basis function are calculated and stored in preprocessing stage. 
\begin{algorithm}[H]
\caption{Inverse Mapping to Subdomain Functions}\label{alg:cap3}
\begin{algorithmic}
  \For{$i \in$ Non-empty Child Box Indices}
    \State{Initialize $\bar{V}_{i}^{F} = 0$ for the $i$\textsuperscript{th} box}
    \For{$k \in$ Subdomain Functions Box $i$}
        \State{Calculate~\eqref{eqn:17} with $\bar{E}_{i}$, to obtain $\bar{Q}^{\prime}_{k}$.}
        \State{$\bar{V}_{i}^{F} = \bar{V}_{i}^{F} + \bar{\bar{D}}_{k} \bar{Q}^{\prime}_{k}$.}
    \EndFor
  \EndFor
  \State \textbf{return} {$\bar{V}_i^{F}$'s.}
\end{algorithmic}
\end{algorithm}

\subsection{Benefits of the Proposed Approach}
The proposed method introduces a three-stage computation step implemented at every iteration of an iterative solver. We now discuss the benefits brought by this novel method.
\begin{itemize}
\item The proposed method achieves a strong scalability for parallelization as each interaction corresponds to distinct regions of the impedance matrix $\bar{\bar{Z}}$, without the need for communication between nodes during interaction computation. Upon completion of the requisite coefficients for mapping parts, these interactions, which are identical for each and every pair of boxes, become the sole components requiring evaluation. The mapping of the boxes into identical structures also ensures a balanced load distribution among nodes. 
\item The proposed method avoids the LFB problem by either directly evaluating the Dyadic Green's function, or training an ML network based on the direct results of the Green's function that does not suffer from the numerical instabilities of Hankel functions in MLFMA implementations.
\item Upon completion of the training process, which is conducted offline, the proposed method can be used in conjunction with any IEs for scattering from conducting or material objects that fit within the mother box. Should the size increase, only additional translation networks require training. Moreover, fast solution of volume-surface IEs~\cite{surfvolIE} for scattering from composite conducting-dielectric objects may also be feasible. Note that as an alternative, even employing the DGFA directly may be considered.
\item The method can be extended to multiple levels, adopting a hierarchical structure akin to the transition from FMM to MLFMA. A tree structure can be constructed within this hierarchical framework, enabling the evaluation of necessary group-by-group interactions at each level, along with the implementation of mapping stages at successive levels. However, uniform basis function locations should be reconsidered from an efficiency point of view. Although strong scalability may not be possible, load balancing for parallelization is expected to be superior to MLFMA due to structurally identical boxes at each level, though further analyses are required. 
\end{itemize}

\section{ML Implementation}\label{sec3}
\begin{figure}[H] 
\centering
\def\svgwidth{3.0in}
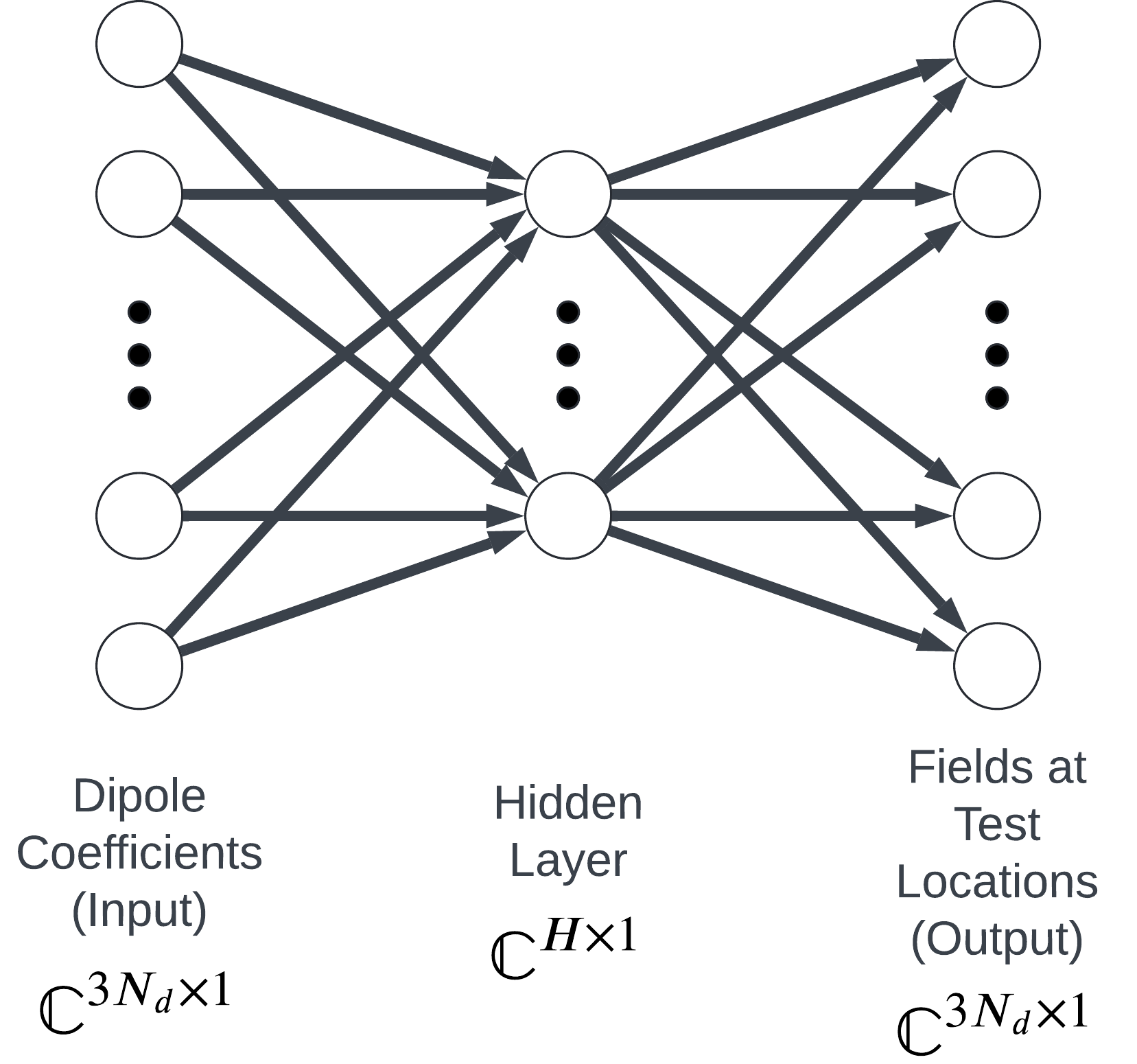
\caption{The structure of the feed-forward and complex-valued fully connected neural network with one hidden layer.}
\label{fig5}
\end{figure}

The primary computational burden of an MVM is concentrated in the translation stage, as defined in Section~\ref{sec2}. Therefore, we employ feed-forward, complex-valued fully connected neural networks (CVFCNNs)~\cite{complex} with strictly lesser complexity than DGFA, which can predict $\bar{\bar{G}}_{\bar{w}_{ij}}$ in~\eqref{eqn:13} accurately. The CVFCNN, shown in Fig.~\ref{fig5}, is constructed with one hidden layer with a dimension of $H$, input and output dimensions of $3N_d \times 1$. This CVFCNN accepts the normalized values of the $3N_d \times 1$ input vector $\bar{P}_{i}$ in~\eqref{eqn:13}, which is normalized with respect to its $\ell_\infty$-norm, and provides the $3N_d \times 1$ output vector $\bar{O}_{i-j}$ that is aimed to be a very accurate prediction of $\bar{E}_{i-j}$ in~\eqref{eqn:13}. The constructed networks are trained using the complex-valued neural networks (CVNN) library built on TensorFlow~\cite{cvnn}, and the Cartesian rectified linear unit (ReLU) activation function is used for the hidden layer. 

An unconventional two-stage training strategy with two different training sets is employed. The training starts with a uniformly generated training dataset. However, because we want to reduce the networks' total training time, the first stage is followed by a second stage with a more realistic dataset that networks may encounter as its inputs during MVMs. That is, the second stage training set is generated as Hertzian dipole coefficients based on pseudo-randomly generated surfaces.

In the first stage, the network undergoes training with a uniformly generated dataset, where real and imaginary parts of all dipole coefficients are randomly assigned to values in the $[-1,1]$ range. This stage continues until the relative $\ell_{2}$-norm error, given by
\begin{equation}
    \label{eqn:171}
    e^{u} = \frac{\| \bar{O}_{i-j}^{u} - \bar{E}_{i-j}\|_{2}}{\| \bar{E}_{i-j}\|_{2}},
\end{equation}
is less than $10^{-3}$. In~\eqref{eqn:171}, $\bar{O}_{i-j}^{u}$ is the output of the test set generated from the uniform distribution, and $\bar{E}_{i-j}$ is the result of~\eqref{eqn:13} computed via DGFA.

\begin{figure}[H] 
\centering
\includegraphics[width=4.5in]{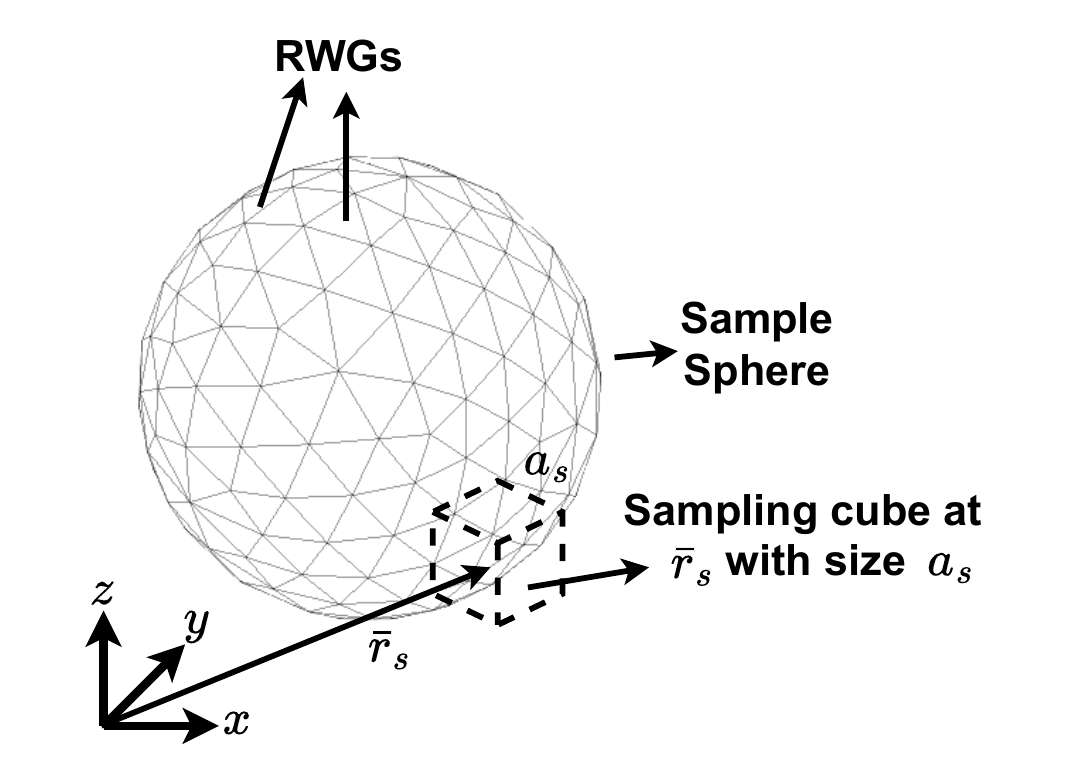}%
\caption{Illustration of pseudo-randomly surface set generation using a PEC sphere discretized with RWGs.}
\label{sphereset}
\end{figure}

\begin{figure}[H]
    \centering
    \begin{subfigure}{2in} 
        \def\svgwidth{\linewidth}
        \relscale{0.6}{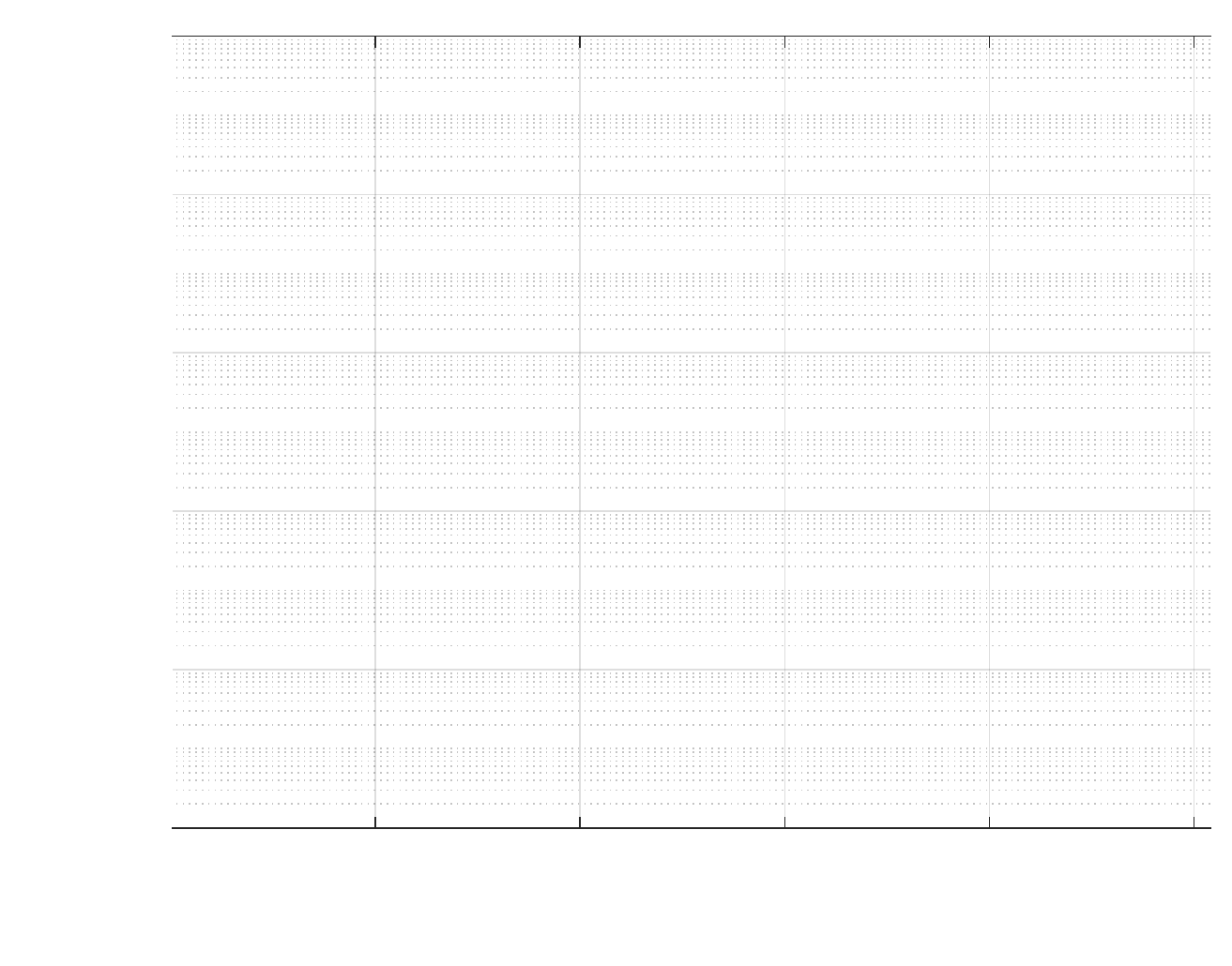} 
        \caption{First stage}
        \label{fig6:sub1}
    \end{subfigure}
    \hspace{0.2\linewidth} 
    \begin{subfigure}{2in} 
        \def\svgwidth{\linewidth}
        \relscale{0.6}{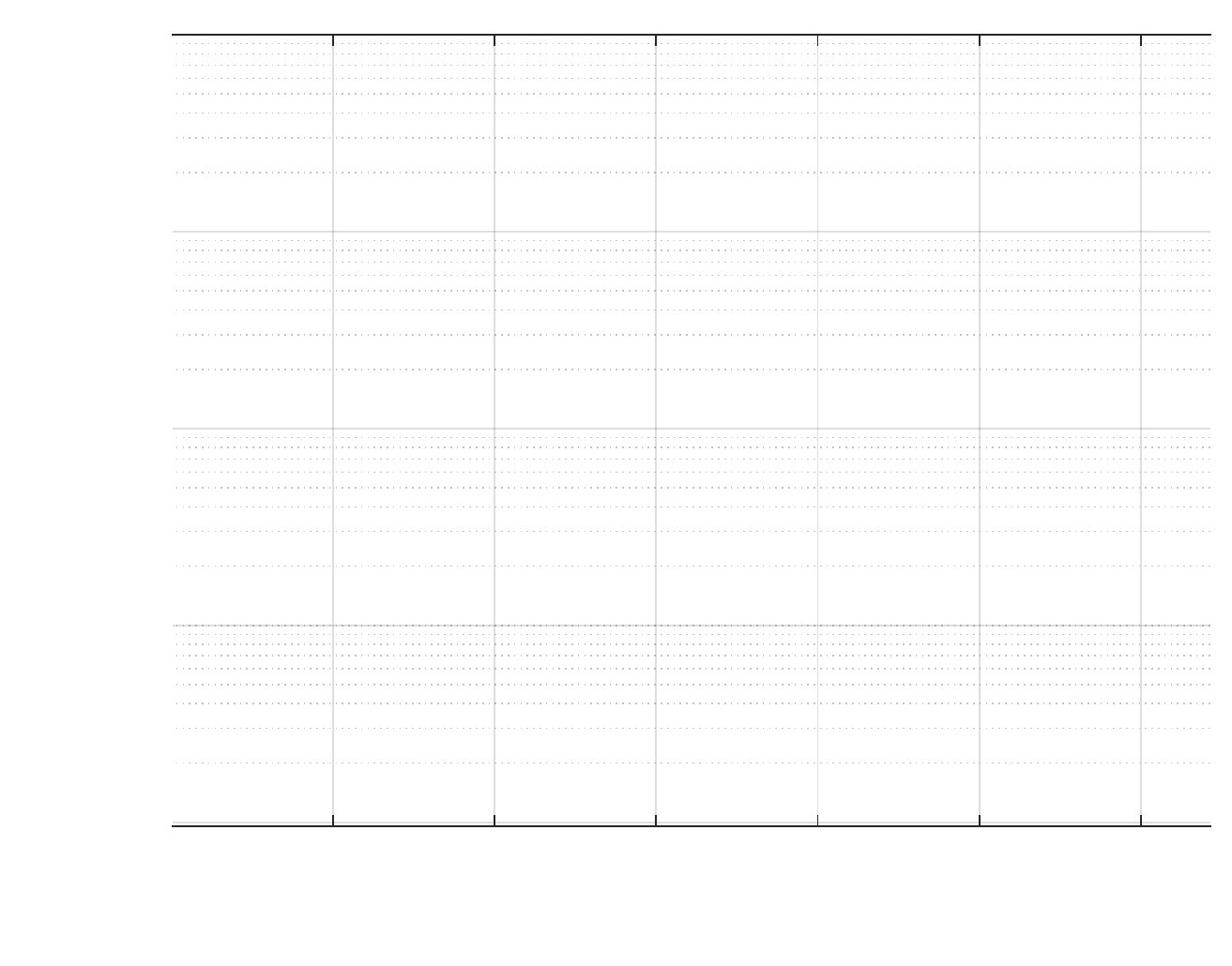} 
        \caption{Second stage}
        \label{fig6:sub2}
    \end{subfigure}
    
    \caption{Training and validation losses at the first and second steps of the training for the $[3\ 3\ 3]a,\ a=\lambda/8$ translation network.}
    \label{fig6}
\end{figure}

In the second stage, the aforementioned pseudo-randomly generated surface dataset is used. The surface generation process is illustrated in~Fig.~\ref{sphereset}. A perfect electric conductor (PEC) sphere discretized with Rao–Wilton-Glisson (RWG) basis functions is sampled to construct training-test sets. Briefly, the center $\bar{r}_{s}$ of a box with length $a_{s}$, which is the size of a child box, is moved across the object with a step size of $0.5a_{s}$ on the surface of the sphere such that the box covers each RWG on the sphere. At each step, complex values randomly selected from the $[-1,1]$ range (for real and imaginary parts) are assigned as the coefficients of the RWG functions that remain within the box. These RWG functions are then mapped to Hertzian dipole coefficients following the aforementioned procedure explained in Section~\ref{sec2}. These Hertzian dipole coefficients form the second training and test set input samples. Similar to the previous stage, the second stage also continues until the relative $\ell_{2}$-norm error, defined as 
\begin{equation}
    \label{eqn:172}
    e^{sp} = \frac{\| \bar{O}_{i-j}^{sp} - \bar{E}_{i-j}\|_{2}}{\| \bar{E}_{i-j}\|_{2}},
\end{equation}
is less than $10^{-3}$. In~\eqref{eqn:172}, $\bar{O}_{i-j}^{sp}$ is the output of the test set generated from the pseudo-randomly generated surface distribution, and $\bar{E}_{i-j}$ is the same as before. 
Throughout the training, the Adam optimizer~\cite{adam} and mean square error loss functions are used for both stages. For demonstration purposes, our CVFCNNs are trained for two different $a$ values: $a=\lambda/8$ to solve any object up to $2\lambda$, and $a=\lambda/256$ to solve any object up to $\lambda/16$, where the conventional MLFMA suffers from the LFB problem. For both cases, hidden layer dimension $H$ is selected as 50, which is a good trade-off between accuracy and low complexity. The input-output dimension $3N_d$ is 192 because we use $4\times 4 \times 4$ dipole grid in each child box, which makes $N_d=64$. Moreover, variable training batch size is used based on experimental outcomes at both stages to achieve lower error rates. The first and second stages of training last approximately 500 and 600 epochs, respectively. Fig.~\ref{fig6} exhibits the training and validation losses of both stages for $a=\lambda/8$ when the translation distance is $[3\ 3\ 3]a$. The loss decreases smoothly in the first stage, but we observe a significantly lower (albeit oscillating) loss in the second stage. This shows that training with only a uniformly generated dataset is not sufficient to reach low error rates and the two-stage training scheme makes the training more efficient and accurate. Moreover, the relative $\ell_2$ errors of several trained networks are listed in Table~\ref{tab1}. Table~\ref{tab1} shows that the $e^{sp}$ value of networks is strictly less than $10^{-3}$, which is a typical relative tolerance selected in iterative solvers. Consequently, the listed networks do not introduce a larger error than the iterative solver tolerance.
\begin{table}[H]
\centering
\caption{Test Set Errors of Several Networks}
\label{tab1}
\begin{tabular}{lcr} 
\toprule
Translation Distance \& & Mean of Relative $\ell_2$ & Std. of Relative $\ell_2$ \\
Input Box Length & Error of Output & Error of Output \\
\midrule
$[0\ 1\ 2]a,\ a=\lambda/8$   & 8.7749e-04   & 4.8696e-04\\ 
$[0\ 0\ 3]a,\ a=\lambda/8$   & 2.0677e-04   & 1.2839e-04      \\
$[3\ 3\ 3]a,\ a=\lambda/8$  & 2.8932e-04   & 2.2387e-04 \\
$[0\ 0\ 2]a,\ a=\lambda/256$   & 4.8532e-04   & 3.0749e-04      \\
$[0\ 2\ 2]a,\ a=\lambda/256$   & 1.1256e-04   & 6.2188e-05      \\
\bottomrule
\end{tabular}
\end{table}

\section{Parallelization}\label{sec4}
As explained in Section~\ref{sec2}, once the mapping to uniform basis functions stage is completed, all boxes become structurally identical, which enables strong scalability for the parallelization process. This is true regardless of our choice to use DGFA or ML networks to calculate the dipole-to-dipole interactions. We now briefly explain the parallelization scheme to achieve a balanced workload among the available workers in a computational platform. Note that all the steps described below can take place in the preprocessing stage. 

Let $\mathbb{W}$ be the set of all translation vectors such that any two translation vectors in this set are not equivalent through rotational and mirror symmetries, and let $\mathbb{W}_{r} \subseteq \mathbb{W}$ denote the set of translation vectors that are not assigned to workers. For each translation vector $\bar{w} \in \mathbb{W}$, let \(n_{\bar{w}}\) be the total number of box-to-box translations with translation vector $\bar{w}$ or its equivalents through rotational and mirror symmetries. Lastly, let $t_{i}$ be the number of translation operations assigned to worker $i$ among a total of $M_{w}$ workers. Then, for the best load balancing as well as to achieve minimum memory allocation for DGFA or ML networks, the following three-step parallelization strategy for worker scheduling is implemented.
\begin{enumerate}
    \item Translation operations belonging to the translation vector with the minimum translation count (i.e., ${\mathrm{argmin}_{\bar{w} \in \mathbb{W}_{r}} \ n_{\bar{w}}}$) are assigned to each worker sequentially in a circular manner. In other words, we assign translations to workers 1 through $M_{w}$, and then start again from the first worker. This stage continues until the number of remaining box-to-box translations, denoted by $r$, is less than the imbalance among the workers given by 
    \begin{equation}
        \label{eqn:+}
        r = \sum_{\bar{w}\in \mathbb{W}_{r}}  n_{\bar{w}}<\sum_{j} \max{\{t_{i} | \, \forall i\}} - t_{j}.
    \end{equation}
    When~\eqref{eqn:+} is fulfilled, this step ends prematurely without adding the last translation vector. Note that there are no memory duplications (e.g., weights for ML networks) in this step, as all equivalent translations are assigned to separate workers.
    \item The remaining translation vectors (i.e., ${\bar{w} \in \mathbb{W}_{r}}$) are concatenated into a list in a way that corresponding ${n}_{\bar{w}}$'s are in ascending order. Then, starting with the worker with the minimum \(t_{j}\), remaining translations are allocated to each worker until the added translation count is below the difference \(\max{\{t_{i} | \, \forall i\}} - t_{j}\), and the same process is applied for all workers once except worker $\mathrm{argmax}_{i}\{{t}_{i} | \, \forall i\}$. This step may induce a slight memory increase due to the orderly assignment of translations.
    \item The residual translations are distributed among all workers to minimize inhomogeneity, with the final stage causing minimal memory increase due to the small size of remaining translations.
\end{enumerate}
Note that in this study, we choose the same ML network structure for all translation vectors. Thus, scheduling can be done only through the translation counts. For optimum implementation of ML networks, the complexity of each network may differ. In this case, the difference in ML network complexities should be considered for scheduling and load balancing. 

Finally, we summarize the MVM process in the iterative solver in Algorithm~\ref{alg:cap5}, where each worker calculates the assigned translations in parallel using ML networks.

\begin{algorithm}[H]
\caption{MVM Part}\label{alg:cap5}
\begin{algorithmic}
\Procedure{MVM Stage}{$x$}
  \State{Obtain uniform basis function values $\bar{P}_i$'s for each box}
  \State{with $x$ and applying Algorithm~\ref{alg:cap1}.}
  \ParFor{$i=1,\ldots,M_{w}$}
    \For{$j=1,\ldots,t_{i}$}
        \State{Calculate $j^{th}$ box-to-box translation in $i^{th}$}
        \State{worker with ML networks, and obtain $\bar{E}_{k}$'s.}
    \EndFor
  \EndParFor
  \State{Calculate test function value in the related box with $\bar{E}_{k}$'s}
  \State{and applying Algorithm~\ref{alg:cap3}.}
\EndProcedure
\end{algorithmic}
\end{algorithm}

\begin{figure}[H]
\centering
\def\svgwidth{3.49in}
\relscale{0.8}{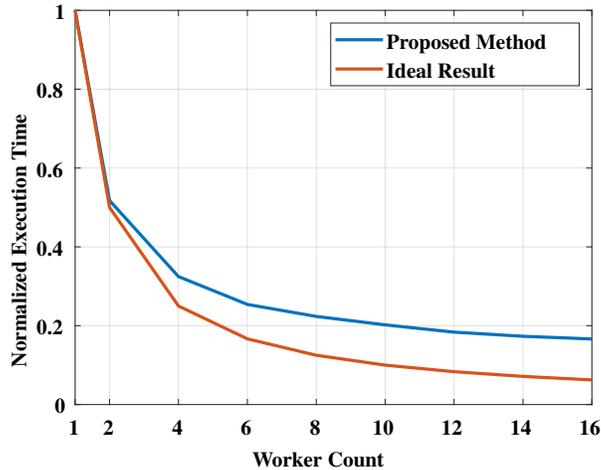} 
\caption{Normalized execution time of the proposed method for the far-zone MVMs of a 0.6m diameter PEC sphere discretized with 114318 basis functions. Illumination is a unit amplitude uniform plane wave at 1 GHz.}
\label{fig7}
\end{figure}

We now demonstrate the strong scalability of the proposed method when ML networks are used for translations. The average MVM computation time for the far-zone interaction calculation of a PEC sphere, discretized with 114318 RWG basis functions, is considered. The PEC sphere is divided into $8^{4}$ boxes so that each child box can have the same size, $a$, with that of the leaf-level boxes when 5-level conventional MLFMA is considered. The normalized execution time versus worker count up to 16 workers is presented and compared with the theoretically ideal situation in Fig.~\ref{fig7} by using Intel Xeon w5-3423 CPU. 
Although the parallelization shows good scalability, the efficiency deviates from the ideal case as the worker count increases. This can be attributed to the fact that our proposed algorithm still has serial computation parts, which are the mapping parts. As the number of workers increases, the execution time of serial computation parts becomes significant, so deviations from the ideal case can arise. Alternatively, the mapping parts can be taken inside of the parallelization process so that strong scaling efficiency increases at the cost of increased memory and time for a small number of workers. Nevertheless, for this example, 5-level MLFMA completes all far-zone interactions in 16.91 seconds, whereas it takes 15.87 seconds for the proposed method when only 4 workers are used. 

\section{Numerical Results} \label{sec5}
Numerical results for PEC and penetrable objects discretized with RWG and Schaubert-Wilton-Glisson (SWG) basis functions, respectively, are provided to assess the performance of the proposed method. The electric field integral equation (D-type formulation for penetrable objects) is used in all examples, where the solutions are obtained using generalized minimal residual (GMRES) iterative solver without any preconditioning with the following parameters: the error tolerance is $10^{-3}$, the restart value is 100, and the maximum number of iterations is fixed to 1000. A second-order Gaussian quadrature is used in the integrations over the surfaces of basis function-related triangles, while a first-order Gaussian quadrature is used in the integrations over the surfaces of testing function-related triangles. Scattering problems solved via the proposed method consist of objects that are PEC spheres, dielectric spheres, a PEC NASA almond model, and a PEC flamme model~\cite{flamme}. All scatterers except the flamme are illuminated by a unit amplitude (1 V/m) $x$-polarized uniform plane wave propagating along the $z$ direction. In the flamme model, a unit amplitude $y$-polarized plane wave impinges from $(\phi,\theta) = (0,60^{\circ})$ that corresponds to $30^{\circ}$ elevation angle from the nose of the flamme. 

\begin{figure}[H]
\centering
\def\svgwidth{4in}
\relscale{0.8}{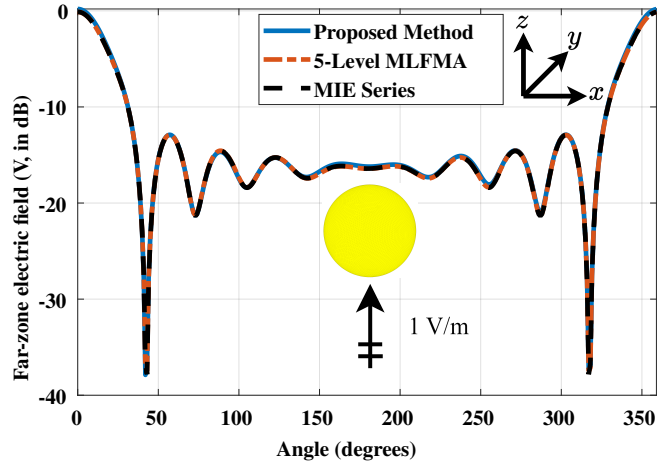} %
\caption{Far-zone co-polarized electric fields scattering from a PEC sphere of diameter 0.6m ($2\lambda$) illuminated with a unit amplitude uniform plane wave at 1 GHz. The discretized surface consists of 114318 basis functions.}
\label{fig8}
\end{figure}

A PEC sphere with a diameter of 0.6m ($2\lambda$) at 1 GHz, discretized with 114318 RWG basis functions, is given as our first example. Note that the same sphere is used in the two-stage training procedure given in Section~\ref{sec3} to train the ML networks. The PEC sphere is divided into $8^{4}$ boxes yielding $a=\lambda/8$ as the size of the child boxes, which is the same as that of the leaf-level boxes when 5-level MLFMA is considered. Fig.~\ref{fig8} presents the far-zone scattered electric field results obtained with the proposed method and is compared with those of 5-level MLFMA and Mie series solutions. The agreement of the proposed method results with the Mie series and 5-level MLFMA solutions is perfect. The proposed method requires 654 iterations to converge, compared to 689 for the 5-level MLFMA. Moreover, one MVM iteration takes 31.55 and 30.49 seconds in MLFMA and our proposed method with 4 worker usage, respectively, using an Intel Xeon w5-3423 CPU.

It should be emphasized at this point that for all the upcoming numerical examples that fit within a mother box of $2\lambda$ and with a child box size of $a=\lambda/8$, the same ML network weights are used. 

\begin{figure}[H]
\centering
\def\svgwidth{4in}
\relscale{0.8}{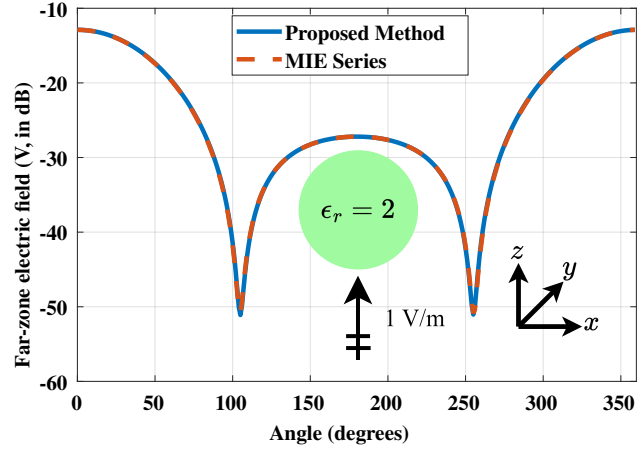} %
\caption{Far-zone co-polarized electric fields scattering from a 0.6m diameter dielectric sphere with $\epsilon_{r} = 2$, illuminated with a unit amplitude uniform plane wave at 250 MHz. The volume of the sphere is discretized with 38169 SWG functions.}
\label{fig9}
\end{figure}

To demonstrate the fact that the trained ML networks can be used regardless of the IE formulation, two penetrable spheres with VIE formulation are considered. First, a dielectric sphere with a diameter of 0.6m ($\lambda/2$) at 250 MHz and $\epsilon_{r}=2$ is discretized with 38169 SWG basis functions and its far-zone scattered electric fields obtained with the proposed method is compared with the Mie series solution to verify the accuracy of the method. As can be seen from Fig.~\ref{fig9}, excellent agreement is observed. 

\begin{figure}[H]
\centering
\def\svgwidth{4in}
\relscale{0.8}{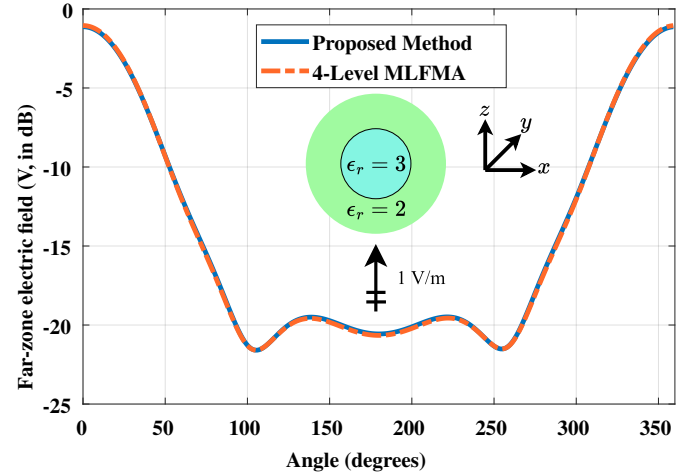} %
\caption{Far-zone co-polarized electric fields scattering from a two-layer dielectric sphere with the core diameter 0.3 m ($\epsilon_{r} = 3$) and the thickness is 0.15 m ($\epsilon_{r} = 2$), illuminated with a unit amplitude uniform plane wave at 500 MHz. The entire volume is discretized with 244689 SWG functions.}
\label{fig10}
\end{figure}

Next, a two-layer dielectric sphere is investigated. Its core region has a diameter of 0.3 m ($\lambda/2$ at 500 MHz) with $\epsilon_{r}=3$, and the outer layer with $\epsilon_{r}=2$ has a thickness of 0.15 m, making the entire scatter's diameter $\lambda$ at 500 MHz (0.6 m). 244689 SWG basis functions are used to discretize this object. The far-zone scattered electric field results obtained with the proposed method as well as with a 4-level MLFMA are compared, as shown in Fig.~\ref{fig10}, where the agreement is excellent again. The required number of iterations for the proposed method and the 4-level MLFMA is 38 and 54, respectively. 

\begin{figure}[H]
\centering
\def\svgwidth{4in}
\relscale{0.8}{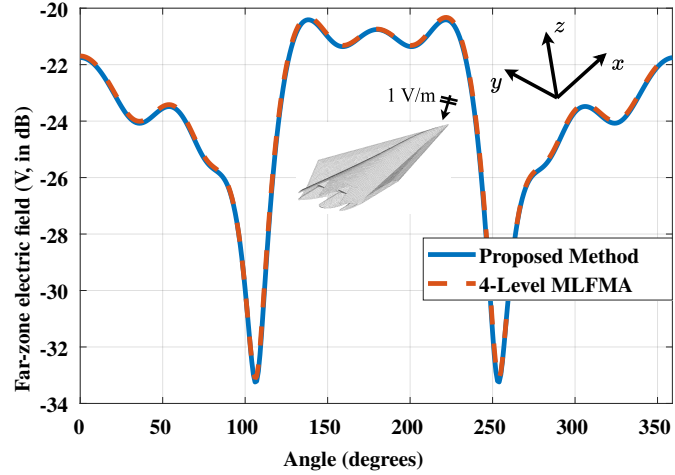} %
\caption{Far-zone electric field along the E-plane ($\phi = 0^{\circ}$) scattered from a 0.6 m PEC flamme model illuminated with a unit amplitude uniform plane wave impinging from $(\phi,\theta) = (0,60^{\circ})$, i.e., $30^{\circ}$ elevation from the nose at 1 GHz. The surface of the flamme is discretized with 51669 RWG functions.}
\label{fig11}
\end{figure}

\begin{figure}[H]
    \centering
    \begin{subfigure}{4.5in} 
        \def\svgwidth{4in}
        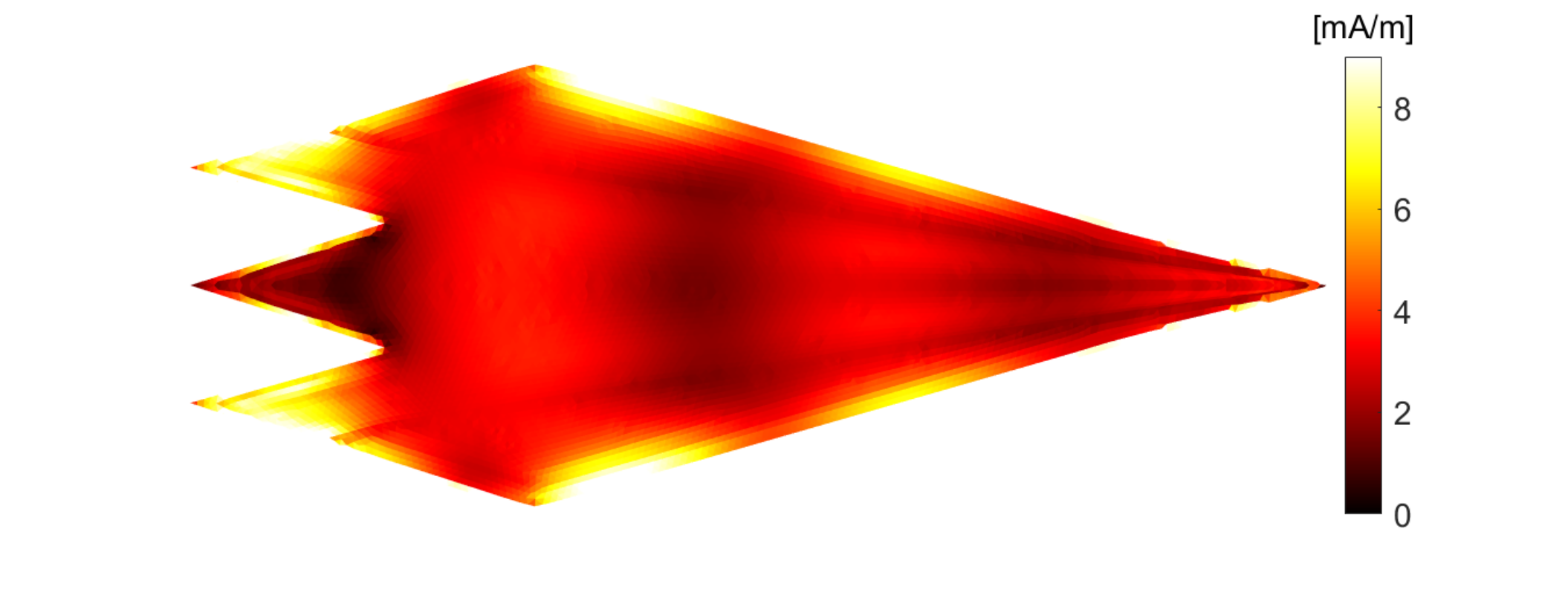 %
        \caption{Proposed Method}
        \label{fig12:sub1}
    \end{subfigure}
    \hfill 
    \begin{subfigure}{4.5in} 
        \def\svgwidth{4in}
        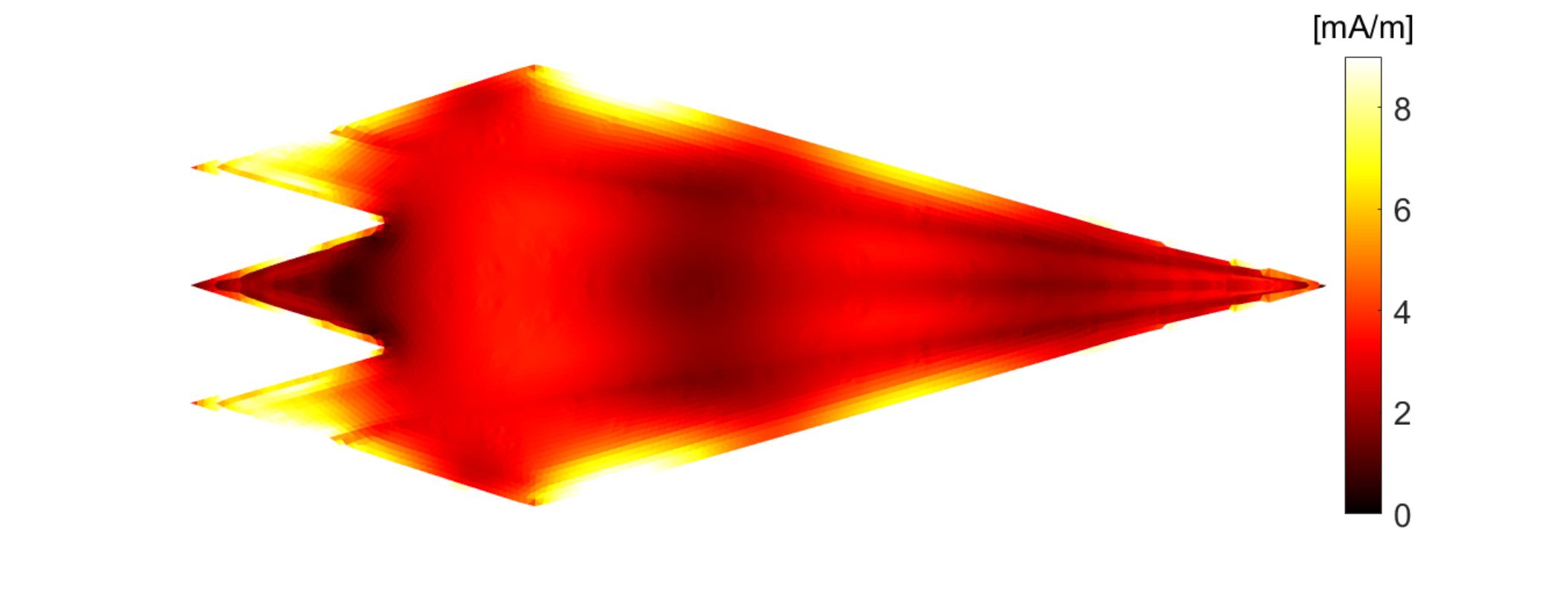 %
        \caption{4-level MLFMA}
        \label{fig12:sub2}
    \end{subfigure}
    
    \caption{Comparison of the induced current density magnitudes for the \textit{(a)} proposed method and \textit{(b)} 4-level MLFMA, when the PEC flamme is illuminated with a unit amplitude uniform plane wave impinging from $(\phi,\theta) = (0,60^{\circ})$, i.e., $30^{\circ}$ elevation from the nose at 1 GHz.}
    \label{fig12}
\end{figure}

The next set of examples are a PEC flamme model and a PEC NASA almond model, both of which have a size of $2\lambda$ at 1 GHz. Note that the same ML network weights from the previous examples are used for both models. The flamme is discretized with 51669 RWG basis functions, while this number is 105159 for the NASA almond. For both scatterers, the E-plane far-zone scattered electric field and the surface current density results are obtained with the proposed method and are compared with those of 4-level MLFMA. Fig.~\ref{fig11} and Fig.~\ref{fig12} show the E-plane far-zone scattered electric field and the surface current density comparisons for the flamme model, respectively, and excellent agreement is observed. Furthermore, one MVM iteration takes 140.45 seconds in our proposed method with 1 worker and 414.6 seconds in 4-level MLFMA by using an Intel Core i7-12700H CPU, with both methods converging within 1200 iterations.

\begin{figure}[H]
\centering
\def\svgwidth{4in}
\relscale{0.8}{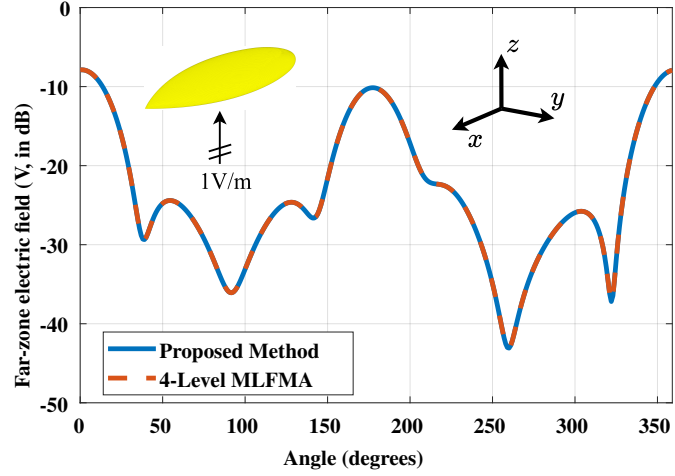} %
\caption{Far-zone electric field along the E-plane scattered from a 0.6 m PEC NASA almond model illuminated with a unit amplitude uniform plane at 1 GHz. The surface of the NASA almond is discretized with 105159 RWG functions.}
\label{fig13}
\end{figure}

\begin{figure}[H]
    \centering
    \begin{subfigure}{4.5in} 
        \def\svgwidth{4in}
        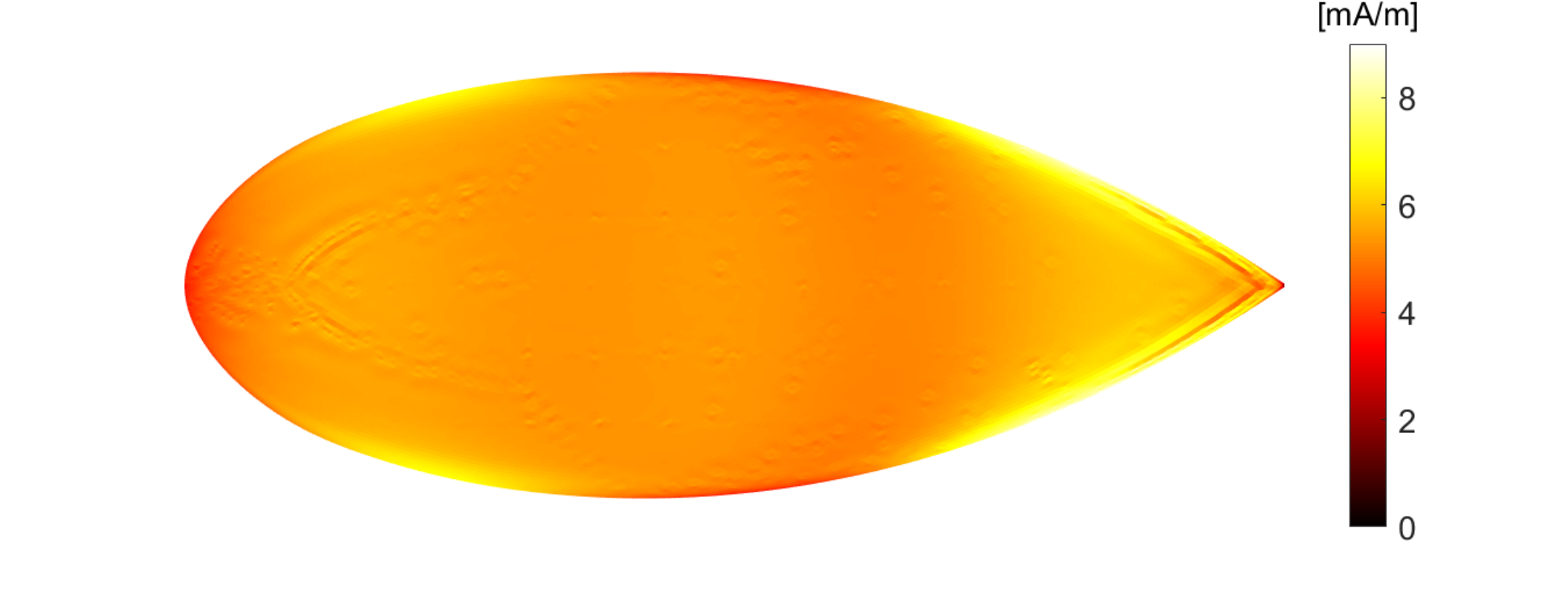 %
        \caption{Proposed Method}
        \label{fig14:sub1}
    \end{subfigure}
    \hfill 
    \begin{subfigure}{4.5in} 
        \def\svgwidth{4in}
        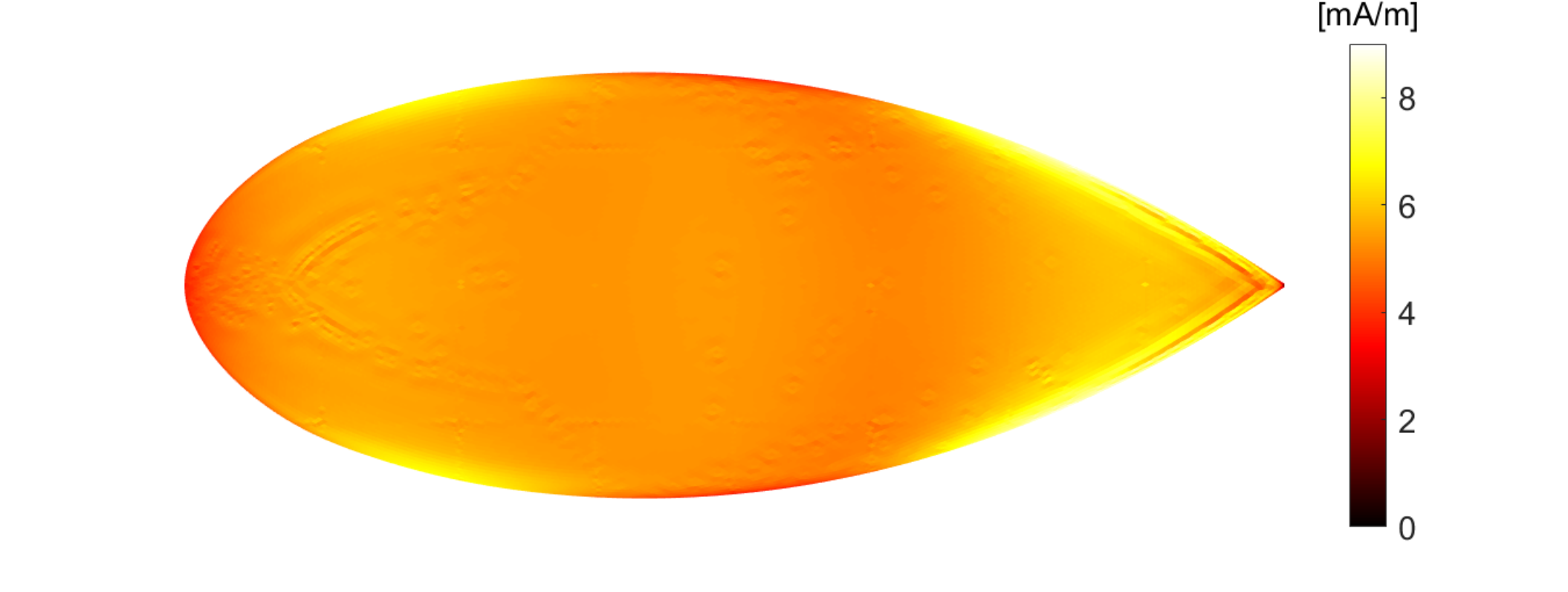 %
        \caption{4-level MLFMA}
        \label{fig14:sub2}
    \end{subfigure}
    
    \caption{Comparison of the induced current density magnitudes for the \textit{(a)} proposed method and \textit{(b)} 4-level MLFMA, when the PEC NASA almond is illuminated with a unit amplitude uniform plane at 1 GHz.}
    \label{fig14}
\end{figure}

Similarly, Fig.~\ref{fig13} and Fig.~\ref{fig14} present the E-plane far-zone scattered electric field and the surface current density comparisons, respectively, for the NASA almond model, and again, very good agreement is observed. Moreover, our proposed method with 1 worker and 4-level MLFMA evaluate one MVM iteration in 423.9 and 1347.3 seconds, respectively, by using Intel Core i7-12700H CPU, with both methods converging in 1900 iterations. 


For the NASA almond and the flamme solutions, our proposed method is faster than 4-level MLFMA even with a single worker. This is because both objects have a high aspect ratio. As a result, the hierarchical structure of MLFMA suffers from inefficiencies, whereas the proposed method can still use symmetries to preserve its efficiency for high aspect ratio objects.

\begin{figure}[H]
\centering
\def\svgwidth{4in}
\relscale{0.8}{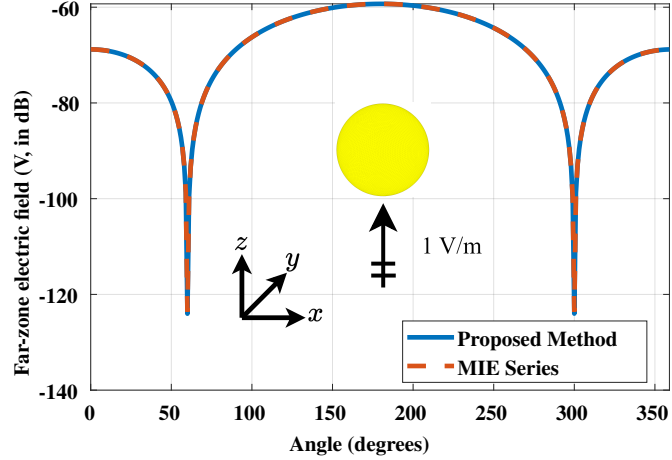} %
\caption{Far-zone co-polarized electric fields scattering from a PEC sphere of diameter 0.6m ($\lambda/64$) illuminated with a unit amplitude uniform plane wave at 7.8125 MHz. The discretized surface consists of 25995 basis functions.}
\label{fig15}
\end{figure}

Finally, our last numerical example demonstrates the proposed method's resilience against the LFB problem that the conventional MLFMA suffers from. For this purpose, we continue to work with the same 0.6 m diameter PEC sphere, but the frequency is set to 7.8125 MHz. At this frequency, the object size becomes $\lambda/64$, and by dividing it to $8^{2}$, the child box size $a$ is $\lambda/256$. ML networks are again trained with the previously explained two-stage training strategy. Fig.~\ref{fig15} presents the far-zone scattered electric field result comparison of the proposed method with Mie series solution and both results agree with each other perfectly. As expected, the conventional MLFMA fails to solve this problem because of LFB. 

It should be kept in mind that the proposed method is not a hierarchical one. Therefore, its efficiency comparison with $N$-level MLFMA when $N\geq 4$ may not be fair if both solutions use a single worker for electrically large objects. However, for all these examples, and many others that are not presented here for brevity, the proposed method is faster than 3-level MLFMA (when a single core is used), and efficiency becomes comparable with $N$-level MLFMA ($N\geq 4$) when a large number of workers is available. We are currently working on a multilevel hierarchical extension of the proposed method. 

\section{Conclusion} \label{sec6}
A novel method for the fast and accurate solution of integral equations that are used for scattering problems is presented. The proposed method employs a novel group-by-group interaction strategy for the evaluation of far-zone interactions by using machine learning networks, showcasing efficiency, accuracy, and strong scalability for parallelization. Moreover, it avoids the low-frequency breakdown problem. Currently, its extension to a multilevel form by adopting a hierarchical structure is underway. Uniform basis function locations, as well as its strong scalability feature, are expected to be affected in its multilevel version. Nevertheless, load balancing is still expected to be superior to MLFMA and/or many similar hierarchical methods available in the literature.

\section*{ACKNOWLEDGEMENT}
Enes Koç was supported by Turkcell as part of 5G and Beyond
Joint Graduate Support Programme coordinated by Information and Communication Technologies Authority.

\bibliographystyle{unsrt}  
\bibliography{main}  






\end{document}

%% file: svg-inkscape/fig5_svg-tex.pdf_tex
\begingroup%
  \makeatletter%
  \providecommand\color[2][]{%
    \errmessage{(Inkscape) Color is used for the text in Inkscape, but the package 'color.sty' is not loaded}%
    \renewcommand\color[2][]{}%
  }%
  \providecommand\transparent[1]{%
    \errmessage{(Inkscape) Transparency is used (non-zero) for the text in Inkscape, but the package 'transparent.sty' is not loaded}%
    \renewcommand\transparent[1]{}%
  }%
  \providecommand\rotatebox[2]{#2}%
  \newcommand*\fsize{\dimexpr\f@size pt\relax}%
  \newcommand*\lineheight[1]{\fontsize{\fsize}{#1\fsize}\selectfont}%
  \ifx\svgwidth\undefined%
    \setlength{\unitlength}{795bp}%
    \ifx\svgscale\undefined%
      \relax%
    \else%
      \setlength{\unitlength}{\unitlength * \real{\svgscale}}%
    \fi%
  \else%
    \setlength{\unitlength}{\svgwidth}%
  \fi%
  \global\let\svgwidth\undefined%
  \global\let\svgscale\undefined%
  \makeatother%
  \begin{picture}(1,0.93584906)%
    \lineheight{1}%
    \setlength\tabcolsep{0pt}%
    \put(0,0){\includegraphics[width=\unitlength,page=1]{fig5_svg-tex.pdf}}%
  \end{picture}%
\endgroup%

%% file: svg-inkscape/firststage_svg-tex.pdf_tex
\begingroup%
  \makeatletter%
  \providecommand\color[2][]{%
    \errmessage{(Inkscape) Color is used for the text in Inkscape, but the package 'color.sty' is not loaded}%
    \renewcommand\color[2][]{}%
  }%
  \providecommand\transparent[1]{%
    \errmessage{(Inkscape) Transparency is used (non-zero) for the text in Inkscape, but the package 'transparent.sty' is not loaded}%
    \renewcommand\transparent[1]{}%
  }%
  \providecommand\rotatebox[2]{#2}%
  \newcommand*\fsize{\dimexpr\f@size pt\relax}%
  \newcommand*\lineheight[1]{\fontsize{\fsize}{#1\fsize}\selectfont}%
  \ifx\svgwidth\undefined%
    \setlength{\unitlength}{630bp}%
    \ifx\svgscale\undefined%
      \relax%
    \else%
      \setlength{\unitlength}{\unitlength * \real{\svgscale}}%
    \fi%
  \else%
    \setlength{\unitlength}{\svgwidth}%
  \fi%
  \global\let\svgwidth\undefined%
  \global\let\svgscale\undefined%
  \makeatother%
  \begin{picture}(1,0.78452381)%
    \lineheight{1}%
    \setlength\tabcolsep{0pt}%
    \put(0,0){\includegraphics[width=\unitlength,page=1]{firststage_svg-tex.pdf}}%
    \put(0.27453196,0.0620436){\makebox(0,0)[lt]{\lineheight{1.25}\smash{\begin{tabular}[t]{l}\textbf{100}\end{tabular}}}}%
    \put(0.44068413,0.0620436){\makebox(0,0)[lt]{\lineheight{1.25}\smash{\begin{tabular}[t]{l}\textbf{200}\end{tabular}}}}%
    \put(0.6068363,0.0620436){\makebox(0,0)[lt]{\lineheight{1.25}\smash{\begin{tabular}[t]{l}\textbf{300}\end{tabular}}}}%
    \put(0.77298834,0.0620436){\makebox(0,0)[lt]{\lineheight{1.25}\smash{\begin{tabular}[t]{l}\textbf{400}\end{tabular}}}}%
    \put(0.93914051,0.0620436){\makebox(0,0)[lt]{\lineheight{1.25}\smash{\begin{tabular}[t]{l}\textbf{500}\end{tabular}}}}%
    \put(0.50018154,0.01028376){\makebox(0,0)[lt]{\lineheight{1.25}\smash{\begin{tabular}[t]{l}\textbf{Epoch}\end{tabular}}}}%
    \put(0,0){\includegraphics[width=\unitlength,page=2]{firststage_svg-tex.pdf}}%
    \put(0.06022257,0.0925819){\makebox(0,0)[lt]{\lineheight{1.25}\smash{\begin{tabular}[t]{l}\textbf{10}\end{tabular}}}}%
    \put(0.10163043,0.11069784){\makebox(0,0)[lt]{\lineheight{1.25}\smash{\begin{tabular}[t]{l}\textbf{-2}\end{tabular}}}}%
    \put(0.06928054,0.22068748){\makebox(0,0)[lt]{\lineheight{1.25}\smash{\begin{tabular}[t]{l}\textbf{10}\end{tabular}}}}%
    \put(0.1106884,0.23880343){\makebox(0,0)[lt]{\lineheight{1.25}\smash{\begin{tabular}[t]{l}\textbf{0}\end{tabular}}}}%
    \put(0.06928054,0.35008707){\makebox(0,0)[lt]{\lineheight{1.25}\smash{\begin{tabular}[t]{l}\textbf{10}\end{tabular}}}}%
    \put(0.1106884,0.36820301){\makebox(0,0)[lt]{\lineheight{1.25}\smash{\begin{tabular}[t]{l}\textbf{2}\end{tabular}}}}%
    \put(0.06928054,0.47819265){\makebox(0,0)[lt]{\lineheight{1.25}\smash{\begin{tabular}[t]{l}\textbf{10}\end{tabular}}}}%
    \put(0.1106884,0.4963086){\makebox(0,0)[lt]{\lineheight{1.25}\smash{\begin{tabular}[t]{l}\textbf{4}\end{tabular}}}}%
    \put(0.06928054,0.60759224){\makebox(0,0)[lt]{\lineheight{1.25}\smash{\begin{tabular}[t]{l}\textbf{10}\end{tabular}}}}%
    \put(0.1106884,0.62570818){\makebox(0,0)[lt]{\lineheight{1.25}\smash{\begin{tabular}[t]{l}\textbf{6}\end{tabular}}}}%
    \put(0.06928054,0.73569782){\makebox(0,0)[lt]{\lineheight{1.25}\smash{\begin{tabular}[t]{l}\textbf{10}\end{tabular}}}}%
    \put(0.1106884,0.75381377){\makebox(0,0)[lt]{\lineheight{1.25}\smash{\begin{tabular}[t]{l}\textbf{8}\end{tabular}}}}%
    \put(0.04391822,0.38696621){\rotatebox{90}{\makebox(0,0)[lt]{\lineheight{1.25}\smash{\begin{tabular}[t]{l}\textbf{Loss}\end{tabular}}}}}%
    \put(0,0){\includegraphics[width=\unitlength,page=3]{firststage_svg-tex.pdf}}%
    \put(0.64819823,0.68773432){\makebox(0,0)[lt]{\lineheight{1.25}\smash{\begin{tabular}[t]{l}\textbf{Training Loss}\end{tabular}}}}%
    \put(0,0){\includegraphics[width=\unitlength,page=4]{firststage_svg-tex.pdf}}%
    \put(0.64819823,0.63355782){\makebox(0,0)[lt]{\lineheight{1.25}\smash{\begin{tabular}[t]{l}\textbf{Validation Loss}\end{tabular}}}}%
    \put(0,0){\includegraphics[width=\unitlength,page=5]{firststage_svg-tex.pdf}}%
  \end{picture}%
\endgroup%

%% file: svg-inkscape/secondstage_svg-tex.pdf_tex
\begingroup%
  \makeatletter%
  \providecommand\color[2][]{%
    \errmessage{(Inkscape) Color is used for the text in Inkscape, but the package 'color.sty' is not loaded}%
    \renewcommand\color[2][]{}%
  }%
  \providecommand\transparent[1]{%
    \errmessage{(Inkscape) Transparency is used (non-zero) for the text in Inkscape, but the package 'transparent.sty' is not loaded}%
    \renewcommand\transparent[1]{}%
  }%
  \providecommand\rotatebox[2]{#2}%
  \newcommand*\fsize{\dimexpr\f@size pt\relax}%
  \newcommand*\lineheight[1]{\fontsize{\fsize}{#1\fsize}\selectfont}%
  \ifx\svgwidth\undefined%
    \setlength{\unitlength}{630bp}%
    \ifx\svgscale\undefined%
      \relax%
    \else%
      \setlength{\unitlength}{\unitlength * \real{\svgscale}}%
    \fi%
  \else%
    \setlength{\unitlength}{\svgwidth}%
  \fi%
  \global\let\svgwidth\undefined%
  \global\let\svgscale\undefined%
  \makeatother%
  \begin{picture}(1,0.78452381)%
    \lineheight{1}%
    \setlength\tabcolsep{0pt}%
    \put(0,0){\includegraphics[width=\unitlength,page=1]{secondstage_svg-tex.pdf}}%
    \put(0.2399428,0.06323407){\makebox(0,0)[lt]{\lineheight{1.25}\smash{\begin{tabular}[t]{l}\textbf{100}\end{tabular}}}}%
    \put(0.37115644,0.06323407){\makebox(0,0)[lt]{\lineheight{1.25}\smash{\begin{tabular}[t]{l}\textbf{200}\end{tabular}}}}%
    \put(0.50237007,0.06323407){\makebox(0,0)[lt]{\lineheight{1.25}\smash{\begin{tabular}[t]{l}\textbf{300}\end{tabular}}}}%
    \put(0.63358358,0.06323407){\makebox(0,0)[lt]{\lineheight{1.25}\smash{\begin{tabular}[t]{l}\textbf{400}\end{tabular}}}}%
    \put(0.76479734,0.06323407){\makebox(0,0)[lt]{\lineheight{1.25}\smash{\begin{tabular}[t]{l}\textbf{500}\end{tabular}}}}%
    \put(0.89601085,0.06323407){\makebox(0,0)[lt]{\lineheight{1.25}\smash{\begin{tabular}[t]{l}\textbf{600}\end{tabular}}}}%
    \put(0.50018154,0.01147424){\makebox(0,0)[lt]{\lineheight{1.25}\smash{\begin{tabular}[t]{l}\textbf{Epoch}\end{tabular}}}}%
    \put(0,0){\includegraphics[width=\unitlength,page=2]{secondstage_svg-tex.pdf}}%
    \put(0.06022257,0.09765436){\makebox(0,0)[lt]{\lineheight{1.25}\smash{\begin{tabular}[t]{l}\textbf{10}\end{tabular}}}}%
    \put(0.10163043,0.1157703){\makebox(0,0)[lt]{\lineheight{1.25}\smash{\begin{tabular}[t]{l}\textbf{-5}\end{tabular}}}}%
    \put(0.06022257,0.25681585){\makebox(0,0)[lt]{\lineheight{1.25}\smash{\begin{tabular}[t]{l}\textbf{10}\end{tabular}}}}%
    \put(0.10163043,0.27493179){\makebox(0,0)[lt]{\lineheight{1.25}\smash{\begin{tabular}[t]{l}\textbf{-4}\end{tabular}}}}%
    \put(0.06022257,0.41727133){\makebox(0,0)[lt]{\lineheight{1.25}\smash{\begin{tabular}[t]{l}\textbf{10}\end{tabular}}}}%
    \put(0.10163043,0.43538727){\makebox(0,0)[lt]{\lineheight{1.25}\smash{\begin{tabular}[t]{l}\textbf{-3}\end{tabular}}}}%
    \put(0.06022257,0.57643282){\makebox(0,0)[lt]{\lineheight{1.25}\smash{\begin{tabular}[t]{l}\textbf{10}\end{tabular}}}}%
    \put(0.10163043,0.59454876){\makebox(0,0)[lt]{\lineheight{1.25}\smash{\begin{tabular}[t]{l}\textbf{-2}\end{tabular}}}}%
    \put(0.06022257,0.7368883){\makebox(0,0)[lt]{\lineheight{1.25}\smash{\begin{tabular}[t]{l}\textbf{10}\end{tabular}}}}%
    \put(0.10163043,0.75500424){\makebox(0,0)[lt]{\lineheight{1.25}\smash{\begin{tabular}[t]{l}\textbf{-1}\end{tabular}}}}%
    \put(0.04391822,0.38815668){\rotatebox{90}{\makebox(0,0)[lt]{\lineheight{1.25}\smash{\begin{tabular}[t]{l}\textbf{Loss}\end{tabular}}}}}%
    \put(0,0){\includegraphics[width=\unitlength,page=3]{secondstage_svg-tex.pdf}}%
    \put(0.64819823,0.68892479){\makebox(0,0)[lt]{\lineheight{1.25}\smash{\begin{tabular}[t]{l}\textbf{Training Loss}\end{tabular}}}}%
    \put(0,0){\includegraphics[width=\unitlength,page=4]{secondstage_svg-tex.pdf}}%
    \put(0.64819823,0.63474829){\makebox(0,0)[lt]{\lineheight{1.25}\smash{\begin{tabular}[t]{l}\textbf{Validation Loss}\end{tabular}}}}%
    \put(0,0){\includegraphics[width=\unitlength,page=5]{secondstage_svg-tex.pdf}}%
  \end{picture}%
\endgroup%

%% file: svg-inkscape/eff_svg-tex.pdf_tex
\begingroup%
  \makeatletter%
  \providecommand\color[2][]{%
    \errmessage{(Inkscape) Color is used for the text in Inkscape, but the package 'color.sty' is not loaded}%
    \renewcommand\color[2][]{}%
  }%
  \providecommand\transparent[1]{%
    \errmessage{(Inkscape) Transparency is used (non-zero) for the text in Inkscape, but the package 'transparent.sty' is not loaded}%
    \renewcommand\transparent[1]{}%
  }%
  \providecommand\rotatebox[2]{#2}%
  \newcommand*\fsize{\dimexpr\f@size pt\relax}%
  \newcommand*\lineheight[1]{\fontsize{\fsize}{#1\fsize}\selectfont}%
  \ifx\svgwidth\undefined%
    \setlength{\unitlength}{630bp}%
    \ifx\svgscale\undefined%
      \relax%
    \else%
      \setlength{\unitlength}{\unitlength * \real{\svgscale}}%
    \fi%
  \else%
    \setlength{\unitlength}{\svgwidth}%
  \fi%
  \global\let\svgwidth\undefined%
  \global\let\svgscale\undefined%
  \makeatother%
  \begin{picture}(1,0.75)%
    \lineheight{1}%
    \setlength\tabcolsep{0pt}%
    \put(0,0){\includegraphics[width=\unitlength,page=1]{eff_svg-tex.pdf}}%
    \put(0.1202381,0.05642857){\makebox(0,0)[lt]{\lineheight{1.25}\smash{\begin{tabular}[t]{l}\textbf{1}\end{tabular}}}}%
    \put(0.17190476,0.05642857){\makebox(0,0)[lt]{\lineheight{1.25}\smash{\begin{tabular}[t]{l}\textbf{2}\end{tabular}}}}%
    \put(0.2752381,0.05642857){\makebox(0,0)[lt]{\lineheight{1.25}\smash{\begin{tabular}[t]{l}\textbf{4}\end{tabular}}}}%
    \put(0.37857143,0.05642857){\makebox(0,0)[lt]{\lineheight{1.25}\smash{\begin{tabular}[t]{l}\textbf{6}\end{tabular}}}}%
    \put(0.48190476,0.05642857){\makebox(0,0)[lt]{\lineheight{1.25}\smash{\begin{tabular}[t]{l}\textbf{8}\end{tabular}}}}%
    \put(0.57571429,0.05642857){\makebox(0,0)[lt]{\lineheight{1.25}\smash{\begin{tabular}[t]{l}\textbf{10}\end{tabular}}}}%
    \put(0.67904762,0.05642857){\makebox(0,0)[lt]{\lineheight{1.25}\smash{\begin{tabular}[t]{l}\textbf{12}\end{tabular}}}}%
    \put(0.78238095,0.05642857){\makebox(0,0)[lt]{\lineheight{1.25}\smash{\begin{tabular}[t]{l}\textbf{14}\end{tabular}}}}%
    \put(0.88571429,0.05642857){\makebox(0,0)[lt]{\lineheight{1.25}\smash{\begin{tabular}[t]{l}\textbf{16}\end{tabular}}}}%
    \put(0.39583369,0.00880952){\makebox(0,0)[lt]{\lineheight{1.25}\smash{\begin{tabular}[t]{l}\textbf{Worker Count}\end{tabular}}}}%
    \put(0,0){\includegraphics[width=\unitlength,page=2]{eff_svg-tex.pdf}}%
    \put(0.09928571,0.08928571){\makebox(0,0)[lt]{\lineheight{1.25}\smash{\begin{tabular}[t]{l}\textbf{0}\end{tabular}}}}%
    \put(0.07190476,0.20761905){\makebox(0,0)[lt]{\lineheight{1.25}\smash{\begin{tabular}[t]{l}\textbf{0.2}\end{tabular}}}}%
    \put(0.07190476,0.32595238){\makebox(0,0)[lt]{\lineheight{1.25}\smash{\begin{tabular}[t]{l}\textbf{0.4}\end{tabular}}}}%
    \put(0.07190476,0.44428571){\makebox(0,0)[lt]{\lineheight{1.25}\smash{\begin{tabular}[t]{l}\textbf{0.6}\end{tabular}}}}%
    \put(0.07190476,0.56261905){\makebox(0,0)[lt]{\lineheight{1.25}\smash{\begin{tabular}[t]{l}\textbf{0.8}\end{tabular}}}}%
    \put(0.09928571,0.68095238){\makebox(0,0)[lt]{\lineheight{1.25}\smash{\begin{tabular}[t]{l}\textbf{1}\end{tabular}}}}%
    \put(0.05642857,0.15892881){\rotatebox{90}{\makebox(0,0)[lt]{\lineheight{1.25}\smash{\begin{tabular}[t]{l}\textbf{Normalized Execution Time}\end{tabular}}}}}%
    \put(0,0){\includegraphics[width=\unitlength,page=3]{eff_svg-tex.pdf}}%
    \put(0.59689,0.63682536){\makebox(0,0)[lt]{\lineheight{1.25}\smash{\begin{tabular}[t]{l}\textbf{Proposed Method}\end{tabular}}}}%
    \put(0,0){\includegraphics[width=\unitlength,page=4]{eff_svg-tex.pdf}}%
    \put(0.59689,0.59293655){\makebox(0,0)[lt]{\lineheight{1.25}\smash{\begin{tabular}[t]{l}\textbf{Ideal Result}\end{tabular}}}}%
    \put(0,0){\includegraphics[width=\unitlength,page=5]{eff_svg-tex.pdf}}%
  \end{picture}%
\endgroup%

%% file: svg-inkscape/2lamsphere1_svg-tex.pdf_tex
\begingroup%
  \makeatletter%
  \providecommand\color[2][]{%
    \errmessage{(Inkscape) Color is used for the text in Inkscape, but the package 'color.sty' is not loaded}%
    \renewcommand\color[2][]{}%
  }%
  \providecommand\transparent[1]{%
    \errmessage{(Inkscape) Transparency is used (non-zero) for the text in Inkscape, but the package 'transparent.sty' is not loaded}%
    \renewcommand\transparent[1]{}%
  }%
  \providecommand\rotatebox[2]{#2}%
  \newcommand*\fsize{\dimexpr\f@size pt\relax}%
  \newcommand*\lineheight[1]{\fontsize{\fsize}{#1\fsize}\selectfont}%
  \ifx\svgwidth\undefined%
    \setlength{\unitlength}{720bp}%
    \ifx\svgscale\undefined%
      \relax%
    \else%
      \setlength{\unitlength}{\unitlength * \real{\svgscale}}%
    \fi%
  \else%
    \setlength{\unitlength}{\svgwidth}%
  \fi%
  \global\let\svgwidth\undefined%
  \global\let\svgscale\undefined%
  \makeatother%
  \begin{picture}(1,0.65625)%
    \lineheight{1}%
    \setlength\tabcolsep{0pt}%
    \put(0,0){\includegraphics[width=\unitlength,page=1]{2lamsphere1_svg-tex.pdf}}%
    \put(0.1185627,0.0638058){\makebox(0,0)[lt]{\lineheight{1.25}\smash{\begin{tabular}[t]{l}\textbf{0}\end{tabular}}}}%
    \put(0.2158851,0.0638058){\makebox(0,0)[lt]{\lineheight{1.25}\smash{\begin{tabular}[t]{l}\textbf{50}\end{tabular}}}}%
    \put(0.31371629,0.0638058){\makebox(0,0)[lt]{\lineheight{1.25}\smash{\begin{tabular}[t]{l}\textbf{100}\end{tabular}}}}%
    \put(0.41918101,0.0638058){\makebox(0,0)[lt]{\lineheight{1.25}\smash{\begin{tabular}[t]{l}\textbf{150}\end{tabular}}}}%
    \put(0.52464574,0.0638058){\makebox(0,0)[lt]{\lineheight{1.25}\smash{\begin{tabular}[t]{l}\textbf{200}\end{tabular}}}}%
    \put(0.63011036,0.0638058){\makebox(0,0)[lt]{\lineheight{1.25}\smash{\begin{tabular}[t]{l}\textbf{250}\end{tabular}}}}%
    \put(0.73557509,0.0638058){\makebox(0,0)[lt]{\lineheight{1.25}\smash{\begin{tabular}[t]{l}\textbf{300}\end{tabular}}}}%
    \put(0.84103981,0.0638058){\makebox(0,0)[lt]{\lineheight{1.25}\smash{\begin{tabular}[t]{l}\textbf{350}\end{tabular}}}}%
    \put(0.38675093,0.02309416){\makebox(0,0)[lt]{\lineheight{1.25}\smash{\begin{tabular}[t]{l}\textbf{Angle (degrees)}\end{tabular}}}}%
    \put(0,0){\includegraphics[width=\unitlength,page=2]{2lamsphere1_svg-tex.pdf}}%
    \put(0.0752048,0.09189683){\makebox(0,0)[lt]{\lineheight{1.25}\smash{\begin{tabular}[t]{l}\textbf{-40}\end{tabular}}}}%
    \put(0.0752048,0.21774689){\makebox(0,0)[lt]{\lineheight{1.25}\smash{\begin{tabular}[t]{l}\textbf{-30}\end{tabular}}}}%
    \put(0.0752048,0.34359695){\makebox(0,0)[lt]{\lineheight{1.25}\smash{\begin{tabular}[t]{l}\textbf{-20}\end{tabular}}}}%
    \put(0.0752048,0.46944701){\makebox(0,0)[lt]{\lineheight{1.25}\smash{\begin{tabular}[t]{l}\textbf{-10}\end{tabular}}}}%
    \put(0.10064958,0.59529707){\makebox(0,0)[lt]{\lineheight{1.25}\smash{\begin{tabular}[t]{l}\textbf{0}\end{tabular}}}}%
    \put(0.06197352,0.12141297){\rotatebox{90}{\makebox(0,0)[lt]{\lineheight{1.25}\smash{\begin{tabular}[t]{l}\textbf{Far-zone electric field (V, in dB)}\end{tabular}}}}}%
    \put(0,0){\includegraphics[width=\unitlength,page=3]{2lamsphere1_svg-tex.pdf}}%
    \put(0.435719,0.56914829){\makebox(0,0)[lt]{\lineheight{1.25}\smash{\begin{tabular}[t]{l}\textbf{Proposed Method}\end{tabular}}}}%
    \put(0,0){\includegraphics[width=\unitlength,page=4]{2lamsphere1_svg-tex.pdf}}%
    \put(0.435719,0.53158254){\makebox(0,0)[lt]{\lineheight{1.25}\smash{\begin{tabular}[t]{l}\textbf{5-Level MLFMA}\end{tabular}}}}%
    \put(0,0){\includegraphics[width=\unitlength,page=5]{2lamsphere1_svg-tex.pdf}}%
    \put(0.435719,0.49401679){\makebox(0,0)[lt]{\lineheight{1.25}\smash{\begin{tabular}[t]{l}\textbf{MIE Series}\end{tabular}}}}%
    \put(0,0){\includegraphics[width=\unitlength,page=6]{2lamsphere1_svg-tex.pdf}}%
  \end{picture}%
\endgroup%

%% file: svg-inkscape/dielecsphere1_svg-tex.pdf_tex
\begingroup%
  \makeatletter%
  \providecommand\color[2][]{%
    \errmessage{(Inkscape) Color is used for the text in Inkscape, but the package 'color.sty' is not loaded}%
    \renewcommand\color[2][]{}%
  }%
  \providecommand\transparent[1]{%
    \errmessage{(Inkscape) Transparency is used (non-zero) for the text in Inkscape, but the package 'transparent.sty' is not loaded}%
    \renewcommand\transparent[1]{}%
  }%
  \providecommand\rotatebox[2]{#2}%
  \newcommand*\fsize{\dimexpr\f@size pt\relax}%
  \newcommand*\lineheight[1]{\fontsize{\fsize}{#1\fsize}\selectfont}%
  \ifx\svgwidth\undefined%
    \setlength{\unitlength}{720bp}%
    \ifx\svgscale\undefined%
      \relax%
    \else%
      \setlength{\unitlength}{\unitlength * \real{\svgscale}}%
    \fi%
  \else%
    \setlength{\unitlength}{\svgwidth}%
  \fi%
  \global\let\svgwidth\undefined%
  \global\let\svgscale\undefined%
  \makeatother%
  \begin{picture}(1,0.65625)%
    \lineheight{1}%
    \setlength\tabcolsep{0pt}%
    \put(0,0){\includegraphics[width=\unitlength,page=1]{dielecsphere1_svg-tex.pdf}}%
    \put(0.1145002,0.08403496){\makebox(0,0)[lt]{\lineheight{1.25}\smash{\begin{tabular}[t]{l}\textbf{0}\end{tabular}}}}%
    \put(0.20821939,0.08403496){\makebox(0,0)[lt]{\lineheight{1.25}\smash{\begin{tabular}[t]{l}\textbf{50}\end{tabular}}}}%
    \put(0.30243021,0.08403496){\makebox(0,0)[lt]{\lineheight{1.25}\smash{\begin{tabular}[t]{l}\textbf{100}\end{tabular}}}}%
    \put(0.40401396,0.08403496){\makebox(0,0)[lt]{\lineheight{1.25}\smash{\begin{tabular}[t]{l}\textbf{150}\end{tabular}}}}%
    \put(0.5055977,0.08403496){\makebox(0,0)[lt]{\lineheight{1.25}\smash{\begin{tabular}[t]{l}\textbf{200}\end{tabular}}}}%
    \put(0.60718154,0.08403496){\makebox(0,0)[lt]{\lineheight{1.25}\smash{\begin{tabular}[t]{l}\textbf{250}\end{tabular}}}}%
    \put(0.70876528,0.08403496){\makebox(0,0)[lt]{\lineheight{1.25}\smash{\begin{tabular}[t]{l}\textbf{300}\end{tabular}}}}%
    \put(0.81034902,0.08403496){\makebox(0,0)[lt]{\lineheight{1.25}\smash{\begin{tabular}[t]{l}\textbf{350}\end{tabular}}}}%
    \put(0.37353912,0.04471221){\makebox(0,0)[lt]{\lineheight{1.25}\smash{\begin{tabular}[t]{l}\textbf{Angle (degrees)}\end{tabular}}}}%
    \put(0,0){\includegraphics[width=\unitlength,page=2]{dielecsphere1_svg-tex.pdf}}%
    \put(0.07262147,0.11116766){\makebox(0,0)[lt]{\lineheight{1.25}\smash{\begin{tabular}[t]{l}\textbf{-60}\end{tabular}}}}%
    \put(0.07262147,0.2088847){\makebox(0,0)[lt]{\lineheight{1.25}\smash{\begin{tabular}[t]{l}\textbf{-50}\end{tabular}}}}%
    \put(0.07262147,0.30660173){\makebox(0,0)[lt]{\lineheight{1.25}\smash{\begin{tabular}[t]{l}\textbf{-40}\end{tabular}}}}%
    \put(0.07262147,0.40431877){\makebox(0,0)[lt]{\lineheight{1.25}\smash{\begin{tabular}[t]{l}\textbf{-30}\end{tabular}}}}%
    \put(0.07262147,0.50203581){\makebox(0,0)[lt]{\lineheight{1.25}\smash{\begin{tabular}[t]{l}\textbf{-20}\end{tabular}}}}%
    \put(0.07262147,0.59975284){\makebox(0,0)[lt]{\lineheight{1.25}\smash{\begin{tabular}[t]{l}\textbf{-10}\end{tabular}}}}%
    \put(0.05984158,0.13967685){\rotatebox{90}{\makebox(0,0)[lt]{\lineheight{1.25}\smash{\begin{tabular}[t]{l}\textbf{Far-zone electric field (V, in dB)}\end{tabular}}}}}%
    \put(0,0){\includegraphics[width=\unitlength,page=3]{dielecsphere1_svg-tex.pdf}}%
    \put(0.39429376,0.56921214){\makebox(0,0)[lt]{\lineheight{1.25}\smash{\begin{tabular}[t]{l}\textbf{Proposed Method}\end{tabular}}}}%
    \put(0,0){\includegraphics[width=\unitlength,page=4]{dielecsphere1_svg-tex.pdf}}%
    \put(0.39429376,0.53296974){\makebox(0,0)[lt]{\lineheight{1.25}\smash{\begin{tabular}[t]{l}\textbf{MIE Series}\end{tabular}}}}%
    \put(0,0){\includegraphics[width=\unitlength,page=5]{dielecsphere1_svg-tex.pdf}}%
  \end{picture}%
\endgroup%

%% file: svg-inkscape/dielecsphere2_svg-tex.pdf_tex
\begingroup%
  \makeatletter%
  \providecommand\color[2][]{%
    \errmessage{(Inkscape) Color is used for the text in Inkscape, but the package 'color.sty' is not loaded}%
    \renewcommand\color[2][]{}%
  }%
  \providecommand\transparent[1]{%
    \errmessage{(Inkscape) Transparency is used (non-zero) for the text in Inkscape, but the package 'transparent.sty' is not loaded}%
    \renewcommand\transparent[1]{}%
  }%
  \providecommand\rotatebox[2]{#2}%
  \newcommand*\fsize{\dimexpr\f@size pt\relax}%
  \newcommand*\lineheight[1]{\fontsize{\fsize}{#1\fsize}\selectfont}%
  \ifx\svgwidth\undefined%
    \setlength{\unitlength}{899.25bp}%
    \ifx\svgscale\undefined%
      \relax%
    \else%
      \setlength{\unitlength}{\unitlength * \real{\svgscale}}%
    \fi%
  \else%
    \setlength{\unitlength}{\svgwidth}%
  \fi%
  \global\let\svgwidth\undefined%
  \global\let\svgscale\undefined%
  \makeatother%
  \begin{picture}(1,0.65638032)%
    \lineheight{1}%
    \setlength\tabcolsep{0pt}%
    \put(0,0){\includegraphics[width=\unitlength,page=1]{dielecsphere2_svg-tex.pdf}}%
    \put(0.12137666,0.05046018){\makebox(0,0)[lt]{\lineheight{1.25}\smash{\begin{tabular}[t]{l}\textbf{0}\end{tabular}}}}%
    \put(0.22061592,0.05046018){\makebox(0,0)[lt]{\lineheight{1.25}\smash{\begin{tabular}[t]{l}\textbf{50}\end{tabular}}}}%
    \put(0.32037577,0.05046018){\makebox(0,0)[lt]{\lineheight{1.25}\smash{\begin{tabular}[t]{l}\textbf{100}\end{tabular}}}}%
    \put(0.4279428,0.05046018){\makebox(0,0)[lt]{\lineheight{1.25}\smash{\begin{tabular}[t]{l}\textbf{150}\end{tabular}}}}%
    \put(0.53550983,0.05046018){\makebox(0,0)[lt]{\lineheight{1.25}\smash{\begin{tabular}[t]{l}\textbf{200}\end{tabular}}}}%
    \put(0.64307696,0.05046018){\makebox(0,0)[lt]{\lineheight{1.25}\smash{\begin{tabular}[t]{l}\textbf{250}\end{tabular}}}}%
    \put(0.750644,0.05046018){\makebox(0,0)[lt]{\lineheight{1.25}\smash{\begin{tabular}[t]{l}\textbf{300}\end{tabular}}}}%
    \put(0.85821103,0.05046018){\makebox(0,0)[lt]{\lineheight{1.25}\smash{\begin{tabular}[t]{l}\textbf{350}\end{tabular}}}}%
    \put(0.39567299,0.00882132){\makebox(0,0)[lt]{\lineheight{1.25}\smash{\begin{tabular}[t]{l}\textbf{Angle (degrees)}\end{tabular}}}}%
    \put(0,0){\includegraphics[width=\unitlength,page=2]{dielecsphere2_svg-tex.pdf}}%
    \put(0.07703127,0.079191){\makebox(0,0)[lt]{\lineheight{1.25}\smash{\begin{tabular}[t]{l}\textbf{-25}\end{tabular}}}}%
    \put(0.07703127,0.18266358){\makebox(0,0)[lt]{\lineheight{1.25}\smash{\begin{tabular}[t]{l}\textbf{-20}\end{tabular}}}}%
    \put(0.07703127,0.28613615){\makebox(0,0)[lt]{\lineheight{1.25}\smash{\begin{tabular}[t]{l}\textbf{-15}\end{tabular}}}}%
    \put(0.07703127,0.38960873){\makebox(0,0)[lt]{\lineheight{1.25}\smash{\begin{tabular}[t]{l}\textbf{-10}\end{tabular}}}}%
    \put(0.09368682,0.49308131){\makebox(0,0)[lt]{\lineheight{1.25}\smash{\begin{tabular}[t]{l}\textbf{-5}\end{tabular}}}}%
    \put(0.10305556,0.59655389){\makebox(0,0)[lt]{\lineheight{1.25}\smash{\begin{tabular}[t]{l}\textbf{0}\end{tabular}}}}%
    \put(0.06349864,0.10937938){\rotatebox{90}{\makebox(0,0)[lt]{\lineheight{1.25}\smash{\begin{tabular}[t]{l}\textbf{Far-zone electric field (V, in dB)}\end{tabular}}}}}%
    \put(0,0){\includegraphics[width=\unitlength,page=3]{dielecsphere2_svg-tex.pdf}}%
    \put(0.43014175,0.55380462){\makebox(0,0)[lt]{\lineheight{1.25}\smash{\begin{tabular}[t]{l}\textbf{Proposed Method}\end{tabular}}}}%
    \put(0,0){\includegraphics[width=\unitlength,page=4]{dielecsphere2_svg-tex.pdf}}%
    \put(0.43014175,0.51542753){\makebox(0,0)[lt]{\lineheight{1.25}\smash{\begin{tabular}[t]{l}\textbf{4-Level MLFMA}\end{tabular}}}}%
    \put(0,0){\includegraphics[width=\unitlength,page=5]{dielecsphere2_svg-tex.pdf}}%
  \end{picture}%
\endgroup%

%% file: svg-inkscape/flamme1_svg-tex.pdf_tex
\begingroup%
  \makeatletter%
  \providecommand\color[2][]{%
    \errmessage{(Inkscape) Color is used for the text in Inkscape, but the package 'color.sty' is not loaded}%
    \renewcommand\color[2][]{}%
  }%
  \providecommand\transparent[1]{%
    \errmessage{(Inkscape) Transparency is used (non-zero) for the text in Inkscape, but the package 'transparent.sty' is not loaded}%
    \renewcommand\transparent[1]{}%
  }%
  \providecommand\rotatebox[2]{#2}%
  \newcommand*\fsize{\dimexpr\f@size pt\relax}%
  \newcommand*\lineheight[1]{\fontsize{\fsize}{#1\fsize}\selectfont}%
  \ifx\svgwidth\undefined%
    \setlength{\unitlength}{1780.5bp}%
    \ifx\svgscale\undefined%
      \relax%
    \else%
      \setlength{\unitlength}{\unitlength * \real{\svgscale}}%
    \fi%
  \else%
    \setlength{\unitlength}{\svgwidth}%
  \fi%
  \global\let\svgwidth\undefined%
  \global\let\svgscale\undefined%
  \makeatother%
  \begin{picture}(1,0.65627633)%
    \lineheight{1}%
    \setlength\tabcolsep{0pt}%
    \put(0,0){\includegraphics[width=\unitlength,page=1]{flamme1_svg-tex.pdf}}%
    \put(0.1216495,0.04968549){\makebox(0,0)[lt]{\lineheight{1.25}\smash{\begin{tabular}[t]{l}\textbf{0}\end{tabular}}}}%
    \put(0.22124274,0.04968549){\makebox(0,0)[lt]{\lineheight{1.25}\smash{\begin{tabular}[t]{l}\textbf{50}\end{tabular}}}}%
    \put(0.32135666,0.04968549){\makebox(0,0)[lt]{\lineheight{1.25}\smash{\begin{tabular}[t]{l}\textbf{100}\end{tabular}}}}%
    \put(0.42928222,0.04968549){\makebox(0,0)[lt]{\lineheight{1.25}\smash{\begin{tabular}[t]{l}\textbf{150}\end{tabular}}}}%
    \put(0.53720778,0.04968549){\makebox(0,0)[lt]{\lineheight{1.25}\smash{\begin{tabular}[t]{l}\textbf{200}\end{tabular}}}}%
    \put(0.64513324,0.04968549){\makebox(0,0)[lt]{\lineheight{1.25}\smash{\begin{tabular}[t]{l}\textbf{250}\end{tabular}}}}%
    \put(0.7530588,0.04968549){\makebox(0,0)[lt]{\lineheight{1.25}\smash{\begin{tabular}[t]{l}\textbf{300}\end{tabular}}}}%
    \put(0.86098437,0.04968549){\makebox(0,0)[lt]{\lineheight{1.25}\smash{\begin{tabular}[t]{l}\textbf{350}\end{tabular}}}}%
    \put(0.39609544,0.00802402){\makebox(0,0)[lt]{\lineheight{1.25}\smash{\begin{tabular}[t]{l}\textbf{Angle (degrees)}\end{tabular}}}}%
    \put(0,0){\includegraphics[width=\unitlength,page=2]{flamme1_svg-tex.pdf}}%
    \put(0.07727992,0.07843198){\makebox(0,0)[lt]{\lineheight{1.25}\smash{\begin{tabular}[t]{l}\textbf{-34}\end{tabular}}}}%
    \put(0.07727992,0.15238127){\makebox(0,0)[lt]{\lineheight{1.25}\smash{\begin{tabular}[t]{l}\textbf{-32}\end{tabular}}}}%
    \put(0.07727992,0.22633057){\makebox(0,0)[lt]{\lineheight{1.25}\smash{\begin{tabular}[t]{l}\textbf{-30}\end{tabular}}}}%
    \put(0.07727992,0.30027987){\makebox(0,0)[lt]{\lineheight{1.25}\smash{\begin{tabular}[t]{l}\textbf{-28}\end{tabular}}}}%
    \put(0.07727992,0.37422917){\makebox(0,0)[lt]{\lineheight{1.25}\smash{\begin{tabular}[t]{l}\textbf{-26}\end{tabular}}}}%
    \put(0.07727992,0.44817846){\makebox(0,0)[lt]{\lineheight{1.25}\smash{\begin{tabular}[t]{l}\textbf{-24}\end{tabular}}}}%
    \put(0.07727992,0.52212776){\makebox(0,0)[lt]{\lineheight{1.25}\smash{\begin{tabular}[t]{l}\textbf{-22}\end{tabular}}}}%
    \put(0.07727992,0.59607706){\makebox(0,0)[lt]{\lineheight{1.25}\smash{\begin{tabular}[t]{l}\textbf{-20}\end{tabular}}}}%
    \put(0.0637399,0.10447067){\rotatebox{90}{\makebox(0,0)[lt]{\lineheight{1.25}\smash{\begin{tabular}[t]{l}\textbf{Far-zone electric field (V,  in dB)}\end{tabular}}}}}%
    \put(0,0){\includegraphics[width=\unitlength,page=3]{flamme1_svg-tex.pdf}}%
    \put(0.65219675,0.27340815){\makebox(0,0)[lt]{\lineheight{1.25}\smash{\begin{tabular}[t]{l}\textbf{Proposed Method}\end{tabular}}}}%
    \put(0,0){\includegraphics[width=\unitlength,page=4]{flamme1_svg-tex.pdf}}%
    \put(0.65219675,0.23549612){\makebox(0,0)[lt]{\lineheight{1.25}\smash{\begin{tabular}[t]{l}\textbf{4-Level MLFMA}\end{tabular}}}}%
    \put(0,0){\includegraphics[width=\unitlength,page=5]{flamme1_svg-tex.pdf}}%
  \end{picture}%
\endgroup%

%% file: svg-inkscape/flammeus_svg-tex.pdf_tex
\begingroup%
  \makeatletter%
  \providecommand\color[2][]{%
    \errmessage{(Inkscape) Color is used for the text in Inkscape, but the package 'color.sty' is not loaded}%
    \renewcommand\color[2][]{}%
  }%
  \providecommand\transparent[1]{%
    \errmessage{(Inkscape) Transparency is used (non-zero) for the text in Inkscape, but the package 'transparent.sty' is not loaded}%
    \renewcommand\transparent[1]{}%
  }%
  \providecommand\rotatebox[2]{#2}%
  \newcommand*\fsize{\dimexpr\f@size pt\relax}%
  \newcommand*\lineheight[1]{\fontsize{\fsize}{#1\fsize}\selectfont}%
  \ifx\svgwidth\undefined%
    \setlength{\unitlength}{1012.5bp}%
    \ifx\svgscale\undefined%
      \relax%
    \else%
      \setlength{\unitlength}{\unitlength * \real{\svgscale}}%
    \fi%
  \else%
    \setlength{\unitlength}{\svgwidth}%
  \fi%
  \global\let\svgwidth\undefined%
  \global\let\svgscale\undefined%
  \makeatother%
  \begin{picture}(1,0.37777778)%
    \lineheight{1}%
    \setlength\tabcolsep{0pt}%
    \put(0,0){\includegraphics[width=\unitlength,page=1]{flammeus_svg-tex.pdf}}%
  \end{picture}%
\endgroup%

%% file: svg-inkscape/flammeMLFMA_svg-tex.pdf_tex
\begingroup%
  \makeatletter%
  \providecommand\color[2][]{%
    \errmessage{(Inkscape) Color is used for the text in Inkscape, but the package 'color.sty' is not loaded}%
    \renewcommand\color[2][]{}%
  }%
  \providecommand\transparent[1]{%
    \errmessage{(Inkscape) Transparency is used (non-zero) for the text in Inkscape, but the package 'transparent.sty' is not loaded}%
    \renewcommand\transparent[1]{}%
  }%
  \providecommand\rotatebox[2]{#2}%
  \newcommand*\fsize{\dimexpr\f@size pt\relax}%
  \newcommand*\lineheight[1]{\fontsize{\fsize}{#1\fsize}\selectfont}%
  \ifx\svgwidth\undefined%
    \setlength{\unitlength}{1012.5bp}%
    \ifx\svgscale\undefined%
      \relax%
    \else%
      \setlength{\unitlength}{\unitlength * \real{\svgscale}}%
    \fi%
  \else%
    \setlength{\unitlength}{\svgwidth}%
  \fi%
  \global\let\svgwidth\undefined%
  \global\let\svgscale\undefined%
  \makeatother%
  \begin{picture}(1,0.37777778)%
    \lineheight{1}%
    \setlength\tabcolsep{0pt}%
    \put(0,0){\includegraphics[width=\unitlength,page=1]{flammeMLFMA_svg-tex.pdf}}%
  \end{picture}%
\endgroup%

%% file: svg-inkscape/badem1_svg-tex.pdf_tex
\begingroup%
  \makeatletter%
  \providecommand\color[2][]{%
    \errmessage{(Inkscape) Color is used for the text in Inkscape, but the package 'color.sty' is not loaded}%
    \renewcommand\color[2][]{}%
  }%
  \providecommand\transparent[1]{%
    \errmessage{(Inkscape) Transparency is used (non-zero) for the text in Inkscape, but the package 'transparent.sty' is not loaded}%
    \renewcommand\transparent[1]{}%
  }%
  \providecommand\rotatebox[2]{#2}%
  \newcommand*\fsize{\dimexpr\f@size pt\relax}%
  \newcommand*\lineheight[1]{\fontsize{\fsize}{#1\fsize}\selectfont}%
  \ifx\svgwidth\undefined%
    \setlength{\unitlength}{1509bp}%
    \ifx\svgscale\undefined%
      \relax%
    \else%
      \setlength{\unitlength}{\unitlength * \real{\svgscale}}%
    \fi%
  \else%
    \setlength{\unitlength}{\svgwidth}%
  \fi%
  \global\let\svgwidth\undefined%
  \global\let\svgscale\undefined%
  \makeatother%
  \begin{picture}(1,0.65606362)%
    \lineheight{1}%
    \setlength\tabcolsep{0pt}%
    \put(0,0){\includegraphics[width=\unitlength,page=1]{badem1_svg-tex.pdf}}%
    \put(0.12159192,0.04960926){\makebox(0,0)[lt]{\lineheight{1.25}\smash{\begin{tabular}[t]{l}\textbf{0}\end{tabular}}}}%
    \put(0.22116909,0.04960926){\makebox(0,0)[lt]{\lineheight{1.25}\smash{\begin{tabular}[t]{l}\textbf{50}\end{tabular}}}}%
    \put(0.32126683,0.04960926){\makebox(0,0)[lt]{\lineheight{1.25}\smash{\begin{tabular}[t]{l}\textbf{100}\end{tabular}}}}%
    \put(0.42917497,0.04960926){\makebox(0,0)[lt]{\lineheight{1.25}\smash{\begin{tabular}[t]{l}\textbf{150}\end{tabular}}}}%
    \put(0.5370831,0.04960926){\makebox(0,0)[lt]{\lineheight{1.25}\smash{\begin{tabular}[t]{l}\textbf{200}\end{tabular}}}}%
    \put(0.64499113,0.04960926){\makebox(0,0)[lt]{\lineheight{1.25}\smash{\begin{tabular}[t]{l}\textbf{250}\end{tabular}}}}%
    \put(0.75289926,0.04960926){\makebox(0,0)[lt]{\lineheight{1.25}\smash{\begin{tabular}[t]{l}\textbf{300}\end{tabular}}}}%
    \put(0.8608074,0.04960926){\makebox(0,0)[lt]{\lineheight{1.25}\smash{\begin{tabular}[t]{l}\textbf{350}\end{tabular}}}}%
    \put(0.39599355,0.00795441){\makebox(0,0)[lt]{\lineheight{1.25}\smash{\begin{tabular}[t]{l}\textbf{Angle (degrees)}\end{tabular}}}}%
    \put(0,0){\includegraphics[width=\unitlength,page=2]{badem1_svg-tex.pdf}}%
    \put(0.07722951,0.07835111){\makebox(0,0)[lt]{\lineheight{1.25}\smash{\begin{tabular}[t]{l}\textbf{-50}\end{tabular}}}}%
    \put(0.07722951,0.18186341){\makebox(0,0)[lt]{\lineheight{1.25}\smash{\begin{tabular}[t]{l}\textbf{-40}\end{tabular}}}}%
    \put(0.07722951,0.2853757){\makebox(0,0)[lt]{\lineheight{1.25}\smash{\begin{tabular}[t]{l}\textbf{-30}\end{tabular}}}}%
    \put(0.07722951,0.388888){\makebox(0,0)[lt]{\lineheight{1.25}\smash{\begin{tabular}[t]{l}\textbf{-20}\end{tabular}}}}%
    \put(0.07722951,0.4924003){\makebox(0,0)[lt]{\lineheight{1.25}\smash{\begin{tabular}[t]{l}\textbf{-10}\end{tabular}}}}%
    \put(0.10326379,0.5959126){\makebox(0,0)[lt]{\lineheight{1.25}\smash{\begin{tabular}[t]{l}\textbf{0}\end{tabular}}}}%
    \put(0.06369168,0.10855108){\rotatebox{90}{\makebox(0,0)[lt]{\lineheight{1.25}\smash{\begin{tabular}[t]{l}\textbf{Far-zone electric field (V, in dB)}\end{tabular}}}}}%
    \put(0,0){\includegraphics[width=\unitlength,page=3]{badem1_svg-tex.pdf}}%
    \put(0.23365638,0.14909489){\makebox(0,0)[lt]{\lineheight{1.25}\smash{\begin{tabular}[t]{l}\textbf{Proposed Method}\end{tabular}}}}%
    \put(0,0){\includegraphics[width=\unitlength,page=4]{badem1_svg-tex.pdf}}%
    \put(0.23365638,0.11070307){\makebox(0,0)[lt]{\lineheight{1.25}\smash{\begin{tabular}[t]{l}\textbf{4-Level MLFMA}\end{tabular}}}}%
    \put(0,0){\includegraphics[width=\unitlength,page=5]{badem1_svg-tex.pdf}}%
  \end{picture}%
\endgroup%

%% file: svg-inkscape/almondus_svg-tex.pdf_tex
\begingroup%
  \makeatletter%
  \providecommand\color[2][]{%
    \errmessage{(Inkscape) Color is used for the text in Inkscape, but the package 'color.sty' is not loaded}%
    \renewcommand\color[2][]{}%
  }%
  \providecommand\transparent[1]{%
    \errmessage{(Inkscape) Transparency is used (non-zero) for the text in Inkscape, but the package 'transparent.sty' is not loaded}%
    \renewcommand\transparent[1]{}%
  }%
  \providecommand\rotatebox[2]{#2}%
  \newcommand*\fsize{\dimexpr\f@size pt\relax}%
  \newcommand*\lineheight[1]{\fontsize{\fsize}{#1\fsize}\selectfont}%
  \ifx\svgwidth\undefined%
    \setlength{\unitlength}{1012.5bp}%
    \ifx\svgscale\undefined%
      \relax%
    \else%
      \setlength{\unitlength}{\unitlength * \real{\svgscale}}%
    \fi%
  \else%
    \setlength{\unitlength}{\svgwidth}%
  \fi%
  \global\let\svgwidth\undefined%
  \global\let\svgscale\undefined%
  \makeatother%
  \begin{picture}(1,0.37777778)%
    \lineheight{1}%
    \setlength\tabcolsep{0pt}%
    \put(0,0){\includegraphics[width=\unitlength,page=1]{almondus_svg-tex.pdf}}%
  \end{picture}%
\endgroup%

%% file: svg-inkscape/almondMLFMA_svg-tex.pdf_tex
\begingroup%
  \makeatletter%
  \providecommand\color[2][]{%
    \errmessage{(Inkscape) Color is used for the text in Inkscape, but the package 'color.sty' is not loaded}%
    \renewcommand\color[2][]{}%
  }%
  \providecommand\transparent[1]{%
    \errmessage{(Inkscape) Transparency is used (non-zero) for the text in Inkscape, but the package 'transparent.sty' is not loaded}%
    \renewcommand\transparent[1]{}%
  }%
  \providecommand\rotatebox[2]{#2}%
  \newcommand*\fsize{\dimexpr\f@size pt\relax}%
  \newcommand*\lineheight[1]{\fontsize{\fsize}{#1\fsize}\selectfont}%
  \ifx\svgwidth\undefined%
    \setlength{\unitlength}{1012.5bp}%
    \ifx\svgscale\undefined%
      \relax%
    \else%
      \setlength{\unitlength}{\unitlength * \real{\svgscale}}%
    \fi%
  \else%
    \setlength{\unitlength}{\svgwidth}%
  \fi%
  \global\let\svgwidth\undefined%
  \global\let\svgscale\undefined%
  \makeatother%
  \begin{picture}(1,0.37777778)%
    \lineheight{1}%
    \setlength\tabcolsep{0pt}%
    \put(0,0){\includegraphics[width=\unitlength,page=1]{almondMLFMA_svg-tex.pdf}}%
  \end{picture}%
\endgroup%

%% file: svg-inkscape/lfbsphere1_svg-tex.pdf_tex
\begingroup%
  \makeatletter%
  \providecommand\color[2][]{%
    \errmessage{(Inkscape) Color is used for the text in Inkscape, but the package 'color.sty' is not loaded}%
    \renewcommand\color[2][]{}%
  }%
  \providecommand\transparent[1]{%
    \errmessage{(Inkscape) Transparency is used (non-zero) for the text in Inkscape, but the package 'transparent.sty' is not loaded}%
    \renewcommand\transparent[1]{}%
  }%
  \providecommand\rotatebox[2]{#2}%
  \newcommand*\fsize{\dimexpr\f@size pt\relax}%
  \newcommand*\lineheight[1]{\fontsize{\fsize}{#1\fsize}\selectfont}%
  \ifx\svgwidth\undefined%
    \setlength{\unitlength}{720bp}%
    \ifx\svgscale\undefined%
      \relax%
    \else%
      \setlength{\unitlength}{\unitlength * \real{\svgscale}}%
    \fi%
  \else%
    \setlength{\unitlength}{\svgwidth}%
  \fi%
  \global\let\svgwidth\undefined%
  \global\let\svgscale\undefined%
  \makeatother%
  \begin{picture}(1,0.65625)%
    \lineheight{1}%
    \setlength\tabcolsep{0pt}%
    \put(0,0){\includegraphics[width=\unitlength,page=1]{lfbsphere1_svg-tex.pdf}}%
    \put(0.1185627,0.0638058){\makebox(0,0)[lt]{\lineheight{1.25}\smash{\begin{tabular}[t]{l}\textbf{0}\end{tabular}}}}%
    \put(0.2158851,0.0638058){\makebox(0,0)[lt]{\lineheight{1.25}\smash{\begin{tabular}[t]{l}\textbf{50}\end{tabular}}}}%
    \put(0.31371629,0.0638058){\makebox(0,0)[lt]{\lineheight{1.25}\smash{\begin{tabular}[t]{l}\textbf{100}\end{tabular}}}}%
    \put(0.41918101,0.0638058){\makebox(0,0)[lt]{\lineheight{1.25}\smash{\begin{tabular}[t]{l}\textbf{150}\end{tabular}}}}%
    \put(0.52464574,0.0638058){\makebox(0,0)[lt]{\lineheight{1.25}\smash{\begin{tabular}[t]{l}\textbf{200}\end{tabular}}}}%
    \put(0.63011036,0.0638058){\makebox(0,0)[lt]{\lineheight{1.25}\smash{\begin{tabular}[t]{l}\textbf{250}\end{tabular}}}}%
    \put(0.73557509,0.0638058){\makebox(0,0)[lt]{\lineheight{1.25}\smash{\begin{tabular}[t]{l}\textbf{300}\end{tabular}}}}%
    \put(0.84103981,0.0638058){\makebox(0,0)[lt]{\lineheight{1.25}\smash{\begin{tabular}[t]{l}\textbf{350}\end{tabular}}}}%
    \put(0.38675093,0.02309416){\makebox(0,0)[lt]{\lineheight{1.25}\smash{\begin{tabular}[t]{l}\textbf{Angle (degrees)}\end{tabular}}}}%
    \put(0,0){\includegraphics[width=\unitlength,page=2]{lfbsphere1_svg-tex.pdf}}%
    \put(0.05892015,0.09189683){\makebox(0,0)[lt]{\lineheight{1.25}\smash{\begin{tabular}[t]{l}\textbf{-140}\end{tabular}}}}%
    \put(0.05892015,0.21725876){\makebox(0,0)[lt]{\lineheight{1.25}\smash{\begin{tabular}[t]{l}\textbf{-120}\end{tabular}}}}%
    \put(0.05892015,0.34262058){\makebox(0,0)[lt]{\lineheight{1.25}\smash{\begin{tabular}[t]{l}\textbf{-100}\end{tabular}}}}%
    \put(0.0752048,0.46798251){\makebox(0,0)[lt]{\lineheight{1.25}\smash{\begin{tabular}[t]{l}\textbf{-80}\end{tabular}}}}%
    \put(0.0752048,0.59334434){\makebox(0,0)[lt]{\lineheight{1.25}\smash{\begin{tabular}[t]{l}\textbf{-60}\end{tabular}}}}%
    \put(0.04568886,0.12141297){\rotatebox{90}{\makebox(0,0)[lt]{\lineheight{1.25}\smash{\begin{tabular}[t]{l}\textbf{Far-zone electric field (V, in dB)}\end{tabular}}}}}%
    \put(0,0){\includegraphics[width=\unitlength,page=3]{lfbsphere1_svg-tex.pdf}}%
    \put(0.60670788,0.16103873){\makebox(0,0)[lt]{\lineheight{1.25}\smash{\begin{tabular}[t]{l}\textbf{Proposed Method}\end{tabular}}}}%
    \put(0,0){\includegraphics[width=\unitlength,page=4]{lfbsphere1_svg-tex.pdf}}%
    \put(0.60670788,0.12351624){\makebox(0,0)[lt]{\lineheight{1.25}\smash{\begin{tabular}[t]{l}\textbf{MIE Series}\end{tabular}}}}%
    \put(0,0){\includegraphics[width=\unitlength,page=5]{lfbsphere1_svg-tex.pdf}}%
  \end{picture}%
\endgroup%

%% file: main.bbl
\begin{thebibliography}{10}

\bibitem{harrington1996field}
Roger~F Harrington and Jan~L Harrington.
\newblock {\em Field computation by moment methods}.
\newblock Oxford University Press, Inc., 1996.

\bibitem{FMM}
R.~Coifman, V.~Rokhlin, and S.~Wandzura.
\newblock The fast multipole method for the wave equation: a pedestrian prescription.
\newblock {\em IEEE Antennas and Propagation Magazine}, 35(3):7--12, 1993.

\bibitem{song1997multilevel}
Jiming Song, Cai-Cheng Lu, and Weng~Cho Chew.
\newblock Multilevel fast multipole algorithm for electromagnetic scattering by large complex objects.
\newblock {\em IEEE transactions on antennas and propagation}, 45(10):1488--1493, 1997.

\bibitem{IL-MLFMA1}
Manouchehr Takrimi, {\"{O}}zg{\"{u}}r Erg{\"{u}}l, and Vakur~B. Ert{\"{u}}rk.
\newblock Incomplete-leaf multilevel fast multipole algorithm for multiscale penetrable objects formulated with volume integral equations.
\newblock {\em IEEE Transactions on Antennas and Propagation}, 65(9):4914--4918, 2017.

\bibitem{IL-MLFMA2}
Bahram Khalichi, {\"{O}}zg{\"{u}}r Erg{\"{u}}l, Manouchehr Takrimi, and Vakur~B. Ert{\"{u}}rk.
\newblock Broadband analysis of multiscale electromagnetic problems: Novel incomplete-leaf {MLFMA} for potential integral equations.
\newblock {\em IEEE Transactions on Antennas and Propagation}, 69(12):9032--9037, 2021.

\bibitem{skeleton}
Yan-Nan Liu, Xiao-Min Pan, and Xin-Qing Sheng.
\newblock Skeletonization accelerated {MLFMA} solution of volume integral equation for plasmonic structures.
\newblock {\em IEEE Transactions on Antennas and Propagation}, 66(3):1590--1594, 2018.

\bibitem{parallel2}
S.~Velamparambil and Weng~Cho Chew.
\newblock Analysis and performance of a distributed memory multilevel fast multipole algorithm.
\newblock {\em IEEE Transactions on Antennas and Propagation}, 53(8):2719--2727, 2005.

\bibitem{parallel}
{\"{O}}zgur Erg{\"u}l and Levent G{\"{u}}rel.
\newblock A hierarchical partitioning strategy for an efficient parallelization of the multilevel fast multipole algorithm.
\newblock {\em IEEE Transactions on Antennas and Propagation}, 57(6):1740--1750, 2009.

\bibitem{Goodfellow}
Ian Goodfellow, Yoshua Bengio, and Aaron Courville.
\newblock {\em Deep Learning}.
\newblock MIT Press, 2016.
\newblock \url{http://www.deeplearningbook.org}.

\bibitem{antenna1}
Qi~Wu, Haiming Wang, and Wei Hong.
\newblock Multistage collaborative machine learning and its application to antenna modeling and optimization.
\newblock {\em IEEE Transactions on Antennas and Propagation}, 68(5):3397--3409, 2020.

\bibitem{antenna2}
Bo~Liu, Hadi Aliakbarian, Zhongkun Ma, Guy A.~E. Vandenbosch, Georges Gielen, and Peter Excell.
\newblock An efficient method for antenna design optimization based on evolutionary computation and machine learning techniques.
\newblock {\em IEEE Transactions on Antennas and Propagation}, 62(1):7--18, 2014.

\bibitem{antenna3}
Jinpil Tak, Adnan Kantemur, Yashika Sharma, and Hao Xin.
\newblock A 3-{D}-printed {W}-band slotted waveguide array antenna optimized using machine learning.
\newblock {\em IEEE Antennas and Wireless Propagation Letters}, 17(11):2008--2012, 2018.

\bibitem{invscat1}
Zhun Wei and Xudong Chen.
\newblock Deep-learning schemes for full-wave nonlinear inverse scattering problems.
\newblock {\em IEEE Transactions on Geoscience and Remote Sensing}, 57(4):1849--1860, 2019.

\bibitem{invscat2}
Huan~Huan Zhang, He~Ming Yao, Lijun Jiang, and Michael Ng.
\newblock Solving electromagnetic inverse scattering problems in inhomogeneous media by deep convolutional encoder–decoder structure.
\newblock {\em IEEE Transactions on Antennas and Propagation}, 71(3):2867--2872, 2023.

\bibitem{intpred}
Yu-Fei Shu, Xing-Chang Wei, Jun Fan, Rui Yang, and Yan-Bin Yang.
\newblock An equivalent dipole model hybrid with artificial neural network for electromagnetic interference prediction.
\newblock {\em IEEE Transactions on Microwave Theory and Techniques}, 67(5):1790--1797, 2019.

\bibitem{time1}
Dal-Jae Yun, In~Il Jung, Haewon Jung, Hoon Kang, Woo-Yong Yang, and In~Yong Park.
\newblock Improvement in computation time of the finite multipole method by using k-means clustering.
\newblock {\em IEEE Antennas and Wireless Propagation Letters}, 18(9):1814--1817, 2019.

\bibitem{time2}
Barışcan Karaosmanoğlu and {\"{O}}zg{\"{u}}r Erg{\"{u}}l.
\newblock Acceleration of {MLFMA} simulations using trimmed tree structures.
\newblock {\em IEEE Transactions on Antennas and Propagation}, 69(1):356--365, 2021.

\bibitem{Halil}
Halil Top{\"{o}}zl{\"{u}}, Barışcan Karaosmanoğlu, and Vakur~Behçet Ert{\"{u}}rk.
\newblock Trimmed multilevel fast multipole algorithm for {D}-type volume integral equations.
\newblock In {\em 2022 IEEE International Symposium on Antennas and Propagation and USNC-URSI Radio Science Meeting (AP-S/URSI)}, pages 1936--1937, 2022.

\bibitem{MLTranslation}
Jia-Jing Sun, Sheng Sun, Yongpin~P. Chen, Lijun Jiang, and Jun Hu.
\newblock Machine-learning-based hybrid method for the multilevel fast multipole algorithm.
\newblock {\em IEEE Antennas and Wireless Propagation Letters}, 19(12):2177--2181, 2020.

\bibitem{LFB1}
L.~Greengard, Jingfang Huang, V.~Rokhlin, and S.~Wandzura.
\newblock Accelerating fast multipole methods for the {Helmholtz} equation at low frequencies.
\newblock {\em IEEE Computational Science and Engineering}, 5(3):32--38, 1998.

\bibitem{LFB2}
J.~S. Zhao and W.~C. Chew.
\newblock Applying {LF–MLFMA} to solve complex {PEC} structures.
\newblock {\em Microwave and Optical Technology Letters}, 28(3):155--160, 2001.

\bibitem{LFB3}
Yun-Hui Chu and Weng~Cho Chew.
\newblock A multilevel fast multipole algorithm for electrically small composite structures.
\newblock {\em Microwave and Optical Technology Letters}, 43(3):202--207, 2004.

\bibitem{LFB5}
Vikram Melapudi, Balasubramaniam Shanker, Sudip Seal, and Srinivas Aluru.
\newblock A scalable parallel wideband {MLFMA} for efficient electromagnetic simulations on large scale clusters.
\newblock {\em IEEE Transactions on Antennas and Propagation}, 59(7):2565--2577, 2011.

\bibitem{LFB6}
L.J. Jiang and W.C. Chew.
\newblock Low-frequency fast inhomogeneous plane-wave algorithm {(LF-FIPWA)}.
\newblock {\em Microwave and Optical Technology Letters}, 40(2):117--122, 2004.

\bibitem{LFB7}
Ignace Bogaert, Joris Peeters, and Femke Olyslager.
\newblock A nondirective plane wave {MLFMA} stable at low frequencies.
\newblock {\em IEEE Transactions on Antennas and Propagation}, 56(12):3752--3767, 2008.

\bibitem{LFB8}
Dennis~T. Schobert and Thomas~F. Eibert.
\newblock Low-frequency surface integral equation solution by multilevel {Green's} function interpolation with fast {Fourier} transform acceleration.
\newblock {\em IEEE Transactions on Antennas and Propagation}, 60(3):1440--1449, 2012.

\bibitem{OzgurBariscan}
{\"{O}}zgur Erg{\"{u}}l and Barışcan Karaosmanoğlu.
\newblock Broadband multilevel fast multipole algorithm based on an approximate diagonalization of the {Green’s} function.
\newblock {\em IEEE Transactions on Antennas and Propagation}, 63(7):3035--3041, 2015.

\bibitem{LFB9}
Mert Kalfa, {\"{O}}zg{\"{u}}r Erg{\"{u}}l, and Vakur~B. Ert{\"{u}}rk.
\newblock Multiple-precision arithmetic implementation of the multilevel fast multipole algorithm.
\newblock {\em IEEE Transactions on Antennas and Propagation}, 72(1):11--21, 2024.

\bibitem{EBF}
S.~Koc, Jiming Song, and W.~C. Chew.
\newblock Error analysis for the numerical evaluation of the diagonal forms of the scalar spherical addition theorem.
\newblock {\em SIAM Journal on Numerical Analysis}, 36(3):906--921, 1999.

\bibitem{EBFLFB}
Mert Kalfa, Vakur~B. Ert{\"{u}}rk, and {\"{O}}zg{\"{u}}r Erg{\"{u}}l.
\newblock Error analysis of {MLFMA} with closed-form expressions.
\newblock {\em IEEE Transactions on Antennas and Propagation}, 69(10):6618--6623, 2021.

\bibitem{Gonz_lez_2009}
\'Alvaro González.
\newblock Measurement of areas on a sphere using {F}ibonacci and latitude–longitude lattices.
\newblock {\em Mathematical Geosciences}, 42(1):49–64, November 2009.

\bibitem{surfvolIE}
Yuan He, Jian~Feng Li, Xiao~Jun Jing, and Mei~Song Tong.
\newblock Fast solution of volume–surface integral equations for multiscale structures.
\newblock {\em IEEE Transactions on Antennas and Propagation}, 67(12):7649--7654, 2019.

\bibitem{complex}
ChiYan Lee, Hideyuki Hasegawa, and Shangce Gao.
\newblock Complex-valued neural networks: A comprehensive survey.
\newblock {\em IEEE/CAA Journal of Automatica Sinica}, 9(8):1406--1426, 2022.

\bibitem{cvnn}
J~Agustin Barrachina.
\newblock Negu93/cvnn: Complex-valued neural networks, November 2022.

\bibitem{adam}
Diederik~P. Kingma and Jimmy Ba.
\newblock Adam: A method for stochastic optimization, 2017.

\bibitem{flamme}
L.~Gürel, H.~Bağcı, J.~C. Castelli, A.~Cheraly, and F.~Tardivel.
\newblock Validation through comparison: Measurement and calculation of the bistatic radar cross section of a stealth target.
\newblock {\em Radio Science}, 38(3), 2003.

\end{thebibliography}
